\documentclass[12pt]{article}
\usepackage{amsmath}
\usepackage{latexsym,amssymb}
\newtheorem{thm}{Theorem}[section]
\newtheorem{lem}[thm]{Lemma}
\newtheorem{prop}[thm]{Proposition}
\newtheorem{cor}[thm]{Corollary}
\newtheorem{df}[thm]{Definition}
\newcommand{\id}{\mathrm{id}}
\newcommand{\Ad}{\mathrm{Ad}\,}
\newcommand{\tM}{\tilde{M}}

\newcommand{\Cnt}{\mathrm{Cnt}}
\newcommand{\Aut}{\mathrm{Aut}}
\newcommand{\Int}{\mathrm{Int}}
\newcommand{\Obm}{\mathrm{Ob}_{\mathrm{m}}}
\newcommand{\md}{\mathrm{mod}}
\newcommand{\tal}{\tilde{\alpha}}
\newcommand{\tbe}{\tilde{\beta}}
\newcommand{\tilu}{\tilde{u}}
\newcommand{\tilU}{\tilde{U}}
\newcommand{\tilG}{\tilde{G}}
\newcommand{\checku}{b}
\newcommand{\cheu}{u}
\newcommand{\baru}{a}
\newcommand{\sr}{\tilde{r}}
\newcommand{\secp}{\tilde{p}}
\newcommand{\sq}{\tilde{q}}
\newcommand{\delv}{\delta v}

\newcommand{\secpq}{\widetilde{pq}}
\newcommand{\secpqr}{\widetilde{pqr}}
\newcommand{\sqr}{\widetilde{qr}}
\newcommand{\circs}{{\lower-.3ex\hbox{{$\scriptscriptstyle\circ$}}}\hspace{1pt}}
\newcommand{\Mopp}{M^{\mathrm{opp}}}

\renewcommand{\epsilon}{\varepsilon}

\author{Toshihiko MASUDA \footnote{Supported by the Grand-in-Aid for Young
Scientists (B), JSPS.} \\
 Graduate School of Mathematics, Kyushu University \\
 744, Motooka, Nishi-ku, 
Fukuoka, 819-0395, JAPAN}
\title{Unified approach to classification of 
 actions of \\ discrete 
amenable groups  on injective
factors}
\date{}

\begin{document}
\maketitle

\begin{abstract}
We 
present a simple unified proof of classification of
 discrete amenable group actions on injective factors.
Our
 argument does not depend on 
 types of factors. We also show the second cohomology vanishing theorem
 for arbitrary cocycle crossed actions of discrete amenable groups 
on the injective factor of type II$_1$.
\end{abstract}

\section{Introduction}\label{sec:intro}
In the theory of operator algebras, the study of automorphism groups and
group actions is one of the most important subjects. 
Especially, 
since Connes' classification of   automorphisms of
the injective factor of type II$_1$ up to
outer conjugacy in  \cite{Co-peri}
and \cite{Con-auto},  
classification of discrete amenable group actions on
injective factors has been developed by many hands.
Namely, in the type II case, V. F. R. Jones classified finite group
actions in \cite{J-act}, and A. Ocneanu developed the method of Connes
and Jones, and classified general discrete
amenable group actions in \cite{Ocn-act}. 
In the type III$_\lambda$ ($\lambda\ne 1$) case, 
C. E. Sutherland and M. Takesaki 
succeed in classification in \cite{Su-Tak-RIMS}, \cite{Su-Tak-act}.
In  the type III$_1$ case,  classification was obtained by 
Y. Kawahigashi,
Sutherland and Takesaki 
for finite
or abelian group actions \cite{KwST}, and finally by Y. Katayama, Sutherland and
Takesaki 
for arbitrary discrete amenable group actions \cite{KtST}.

Although the classification theorem can be stated in unified way as in
\cite{KtST}, proofs presented in the above cited papers heavily
depend on type of 
factors. 
Indeed, in the type III cases, they reduce the problem to that of
the type II case by means of the structure theorem of type III factors. 
(In this approach, type III$_1$ case is extremely hard because of lack
of classification of $\mathbb{R}$-actions. This the reason why type
III$_\lambda $ ($\lambda \ne 1$) and  
III$_1$ case are treated separately.)

In \cite{M-inter-mathscan}, by developing the Evans-Kishimoto intertwining
argument \cite{EvKi}, we presented the classification 
of  centrally free actions of discrete
amenable groups on injective factors, whose  proof  is independent from
types of factors. 
So one may expect if he can give a unified
approach for classification of general actions by a similar technique. 
Indeed, for given two actions of a discrete amenable group $G$ with same
invariants,  if $G$ is a semidirect product $H\rtimes K$, where $H$ is a
centrally trivial part, then
we can show the cocycle
conjugacy of two actions by using \cite{M-inter-mathscan}. 
(See the beginning of \S \ref{sec:cocycle} for details.)

In this paper, we extend the
intertwining argument and present a unified proof of  classification of 
all actions of discrete amenable
groups on injective factors. 
Roughly speaking, we first divide a given action to
a centrally trivial part and a centrally free part. 
It is rather easy to classify centrally trivial parts.
By applying the intertwining argument, we classify centrally free parts,
 and combine  centrally free parts and centrally trivial parts. 
(Similar idea has already appeared in \cite{JT} and \cite{Su-Tak-RIMS}
for classification of groupoid actions.) 
One difficulty is that 
a centrally free part does not give
an action in the usual sense, and this makes our argument more difficult. 
So if we take this fact into account, it is more convenient  
to handle cocycle crossed actions. 
Thus in  the intertwining argument, we also need the second cohomology
vanishing procedure. 
Therefore our argument looks like the  
mixture of  the second cohomology vanishing argument presented in \cite{M-T-CMP}
and  the
intertwining argument. Here we emphasize that our argument is 
based on the Rohlin type theorem \cite{Ocn-act} and the
characterization of approximately inner automorphisms and
centrally trivial automorphisms \cite{Con-surv}, \cite{KwST}. 

This paper is organized as follows. In \S \ref{sec:pre}, we 
recall the definition of invariants of actions and state
the main theorem and its corollaries. In \S \ref{sec:cocycle}
we introduce the notion of quasi cocycle crossed actions. 
In \S \ref{sec:ultra}, by using ultraproduct, we 
show that we can approximate a quasi cocycle crossed action by a
unitary perturbation of another quasi cocycle crossed action.
In \S\ref{sec:coh}, we show the approximate cohomology vanishing, which
is the main tool for intertwining argument. The key ingredient is the
Rohlin type theorem.
In \S \ref{sec:inter}  we classify quasi cocycle crossed actions, and 
show the main theorem by applying the extended intertwining argument.
In \S \ref{sec:model}, we discuss about construction of model actions.
In \S \ref{sec:subfactor}, we treat group actions on subfactors of type
II$_1$ with finite index. By a suitable modification, we can apply the
extended intertwining argument in this case. 
In \S \ref{sec:outer}, we apply the intertwining argument to classify
outer actions in the sense of \cite{KtT-outerI}. 
In appendix, we present
proofs of the second cohomology vanishing theorem, 
the Rohlin theorem and the existence of extension of automorphisms to a 
twisted crossed product von Neumann algebra
for readers' convenience.

The author is much indebted to Professor Katayama for improvement of our
argument, and express his gratitude. 
 He also thanks Professor Takesaki and 
Professor Ueda for
useful comments on this paper.

\section{Preliminaries and main results}\label{sec:pre}
Our standard references for theory of operator algebras are \cite{Tak-book}.
Throughout this paper, $M$ and $G$ always denotes 
an injective factor, and a discrete amenable group respectively, 
although some of results are valid for general factors and discrete
groups. 

\begin{df}
$(1)$ Let $\alpha$ be a map from $G$ into $\Aut(M)$ with $\alpha_e=\id$
 and $v^\alpha(g,h)\in U(M)$, $g,h\in G$, such that $v^\alpha(e,h)=v^\alpha(g,e)=1$. 
We say 
$(\alpha_g, v^\alpha(g,h))$ is a cocycle crossed action of $G$
 if it satisfies 
$$\alpha_g\circs\alpha_h=\Ad v^\alpha(g,h)\circs \alpha_{gh},\,\,\, 
v^\alpha(g,h) v^\alpha(gh,k)=\alpha_g(v^\alpha(h,k))v^\alpha(g,hk).$$
$(2)$ 
Two cocycle crossed actions $(\alpha,v^\alpha(g,h))$ and
$(\beta,v^\beta(g,h))$ are said to be cocycle conjugate if there exist
a family of unitaries $\{u_g\}_{g\in G}$ and $\theta\in \Aut(M)$ such that 
$$\Ad u_g\circs\alpha_g=\theta\circs \beta_g\circs \theta^{-1},\,\,\,
u_g\alpha_g(u_h)v^\alpha(g,h)u_{g h}^*=\theta(v^\beta(g,h)).$$ 
If we can take $\theta\in \overline{\Int}(M)$, then we say they are
strongly cocycle conjugate. 
\end{df}

We generalize the  invariants of group actions \cite{Su-Tak-act},
\cite{KtST} to  cocycle crossed action case. 
Let $\tM$ be the core for $M$ and 
$(\tM,\mathbb{R}, \theta_t)$ the core covariant system \cite{KtST}.
(It is also called the non-commutative flow of weights in \cite{FT}.) 
Denote by $\tal\in \Aut(\tM)$ the canonical extension of $\alpha\in
\Aut(M)$ \cite{HS}. 
Take a faithful normal weight $\psi$. We can identify
$(\tM, \mathbb{R}, \theta)$ 
with $(M\rtimes_{\sigma^\psi}\mathbb{R},
\widehat{\sigma^\psi}_t)$, and $\tal$ is given by 
$$\tal|_M=\alpha,\,\,\, \tal(\lambda_t^\psi)=(D\psi\circs \alpha^{-1}:D\psi)_t\lambda_t^\psi$$
via this identification, where $\lambda^\psi_t$ is the implementing
unitary for $\sigma^\psi_t$.  (See \cite{FT} on the functorial property
of $(\tM,\mathbb{R}, \theta_t)$.)
Let $\tilU(M)$ be the normalizer unitary group for $M\subset \tM$.
Set $$\Cnt_r(M)=\{\sigma\in \Aut(M)\mid \tilde{\sigma} \in \Int(\tM)\}=
\{\Ad u|_M\mid  u\in \tilU(M)\}.$$
Then the normal subgroup of $G$ defined by 
$H=\{g\in G\mid \alpha_g\in \Cnt_r(M)\}$ is the first invariant.
The second invariant is a Connes-Takesaki module
$\mathrm{mod}(\alpha_g)=\tal_g|_{Z(\tM)}$ \cite{CT}.

The third invariant is a characteristic invariant. 
Set $\delv^\alpha(g,n):=v^\alpha(g,g^{-1}ng)v^\alpha(n,g)^*$ for a
2-cocycle $v^\alpha(g,h)$.
It is easy to see $\alpha_g\circs \alpha_{g^{-1}ng}=\Ad
\delv^\alpha(g,n)\,\circs \alpha_n\circs \alpha_g$

Take $\tilde{u}_m^\alpha\in \tilU(M)$ such that $\tilde{\alpha}_m=\Ad
\tilde{u}_m^\alpha$, $m\in H$, and $\tilu_e=1$.  
Comparing 
$$
\tal_g\circs\tal_{g^{-1}n g}\circs\tal_{g}^{-1}=\Ad
\delv^\alpha(g,n)\circs \tal_{n}\circs \tal_g\circs \tal_g^{-1}=
\Ad (\delv^\alpha(g,n)\tilu_n^\alpha )
$$
with 
$$\tal_g\circs\tal_{g^{-1}n g}\circs\tal_{g}^{-1}= \Ad\tal_g(\tilu_{g^{-1}n g}^\alpha),$$
we get a unitary $\lambda(g,n)\in Z(\tM)$ such that
$$\tal_g(\tilu_{g^{-1}n g}^\alpha)=\lambda(g,n)
\delv^\alpha(g,n) \tilu_n^\alpha.$$
In a similar way, from $\tal_{m}\tal_n=\Ad v^\alpha(m,n)\tal_{m n}$, we get 
a unitary $\mu(m,n)\in Z(\tM)$ by
$$\tilde{u}_m^\alpha\tilde{u}_n^\alpha=\mu(m,n)v^\alpha(m,n)\tilde{u}_{m
n}^\alpha,\, m,n\in H.$$

Since $\theta_t$ commutes with $\tilde{\alpha}_g$, we get 
a 1-cocycle $c_t(m)\in Z^1_\theta(\mathbb{R},U(Z(\tM)))$ 
by $\theta_t(\tilde{u}_m^\alpha)=c_{t}(m)\tilde{u}_m^\alpha $. 
Unitaries $(\lambda,\mu)$ and $c_t(n)$ are related as below.
\begin{eqnarray*}
&& \theta_t(\lambda(g,n))c_t(n)=\md(\alpha_g)(c_t(g^{-1}ng))\lambda(g,n), \\
&& c_t(m)c_t(n)\mu(m,n)=c_t(mn)\theta_t(\mu(m,n)).
\end{eqnarray*}
Set $\lambda((g,t),n)=\theta_t(\lambda(g,n))c_t(n)$. 
We extend a cocycle crossed
action $(\tal, v^\alpha(g,h))$ of $G$ to that of $\tilG=G\times \mathbb{R}$ by 
setting $\tal_{(g,t)}=\tal_{g}\theta_t$ and $v^\alpha((g,t),(h,s))=v^\alpha(g,h)$. 
We can verify that  $(\lambda,\mu)$ 
satisfies the following relations. (These include the above relations.)
\begin{eqnarray*}
&& \mu(l,m)\mu(lm,n)=\mu(m,n)\mu(l,mn), \,\, l,m,n\in H, \\
&&\lambda(gh,n)=\tal_g(\lambda(h,g^{-1}ng))\lambda(g,n),\,\, 
g,h\in \tilG, n\in H, \\ 
&&\lambda(g,mn)\lambda(g,m)^*\lambda(g,n)^*=
\mu(m,n)\tal_g(\mu(g^{-1}mg,g^{-1}ng)^*), \,\, g\in \tilG, m,n\in H,\\ 
&&\lambda(m,n)=\mu(m,m^{-1}nm)\mu(n,m)^*,\,\, m,n\in H, \\
&& \lambda(g,n)=\mu(l,m)=1 \mbox{ if any of }g,l,m,n \mbox{ is }e. 
\end{eqnarray*} 
\begin{df}\label{df:chara}
$(1)$ We say  $(\lambda,\mu)$ satisfying the above conditions as a
characteristic cocycle for $\alpha $ and denote the set of all
 characteristic cocycles by   $Z_\alpha(\tilG,H,U(Z(\tM)))$. Note that 
$Z_\alpha(G,H,U(Z(\tM)))$ becomes an abelian group by a natural
multiplication.  \\
$(2)$ Let $B_\alpha(\tilG,H,U(Z(\tM)))$ be a set of coboundaries defined as
$$B_\alpha(\tilG,H,U(Z(\tM)))=\{
(\tal_g(z_{g^{-1}ng})z_n^*, z_mz_nz_{mn}^*)\mid \{z_n\}_{n\in H}\subset
U(Z(\tM)),\,\, z_e=1\},$$
which is a normal subgroup of 
$Z_\alpha(\tilG,H,U(Z(\tM)))$. \\
$(3)$ Let $\Lambda_\alpha(\tilG,H,U(Z(\tM)))$ be a quotient group.  
 The equivalence class $\chi(\alpha)=[\lambda,\nu]\in \Lambda_\alpha(\tilG,H,U(Z(\tM)))$ 
is said to be the characteristic invariant for $\alpha$.
\end{df}
Though $(\lambda,\mu)$ depends on the choice of $\tilu_n$,
$\chi(\alpha)$ does not depend on $\tilu_n$. 
Thus the true invariant for $\alpha$ is $\chi(\alpha)$ rather than
$(\lambda,\mu)$.

We denote the triplet $\{H, \mathrm{mod}(\alpha_g),\chi(\alpha)\}$ by 
$\mathrm{Inv}(\alpha)$. It is a routine work to show
$\mathrm{Inv}(\alpha)$ is a strong cocycle conjugacy invariant. 
The purpose of this paper is to present a proof  the converse
implication which does not depend on types of factors 
by using
the Evans-Kishimoto type intertwining argument as in \cite{M-inter-mathscan}.
\begin{thm}\label{thm:class}
Let $M$ be an injective factor, and $G$ a discrete amenable group.
 Let $(\alpha,v^\alpha(g,h))$ and $(\beta,v^\beta(g,h))$ 
be cocycle crossed actions of $G$ on $M$ with 
 $\mathrm{Inv}(\alpha)=\mathrm{Inv}(\beta)$. 
Then $\alpha$ and $\beta$ are strongly cocycle conjugate.
\end{thm}
The proof of Theorem \ref{thm:class} will be presented in sequel
sections. Here we state corollaries of the main theorem. 

Let $\Aut_{\theta}(Z(\tM))=\{\sigma\in \Aut(Z(\tM))\mid
\sigma\theta_t=\theta_t\sigma\}$. Then $\sigma \in
\Aut_{\theta}(Z(\tM))$ acts on $\mathrm{Inv}(\alpha)$ by 
$\sigma(\mathrm{Inv}(\alpha))=\{H, \sigma\circs\md(\alpha_g)\circs\sigma^{-1},
[\sigma(\lambda),\sigma(\mu)]\}$.
\begin{cor}
Let $M$, $G$ be as in Theorem \ref{thm:class}. 
Let $\alpha$ and $\beta$ be actions of $G$. 
Then $\alpha$ and $\beta $ 
 are cocycle conjugate if and only if there exists $\sigma \in \Aut_{\theta}(Z(\tM))$ 
such that $\mathrm{Inv}(\alpha)=\sigma(\mathrm{Inv}(\beta))$.  
\end{cor}
\textbf{Proof.} This follows from Theorem \ref{thm:class}, 
$\mathrm{Inv}(\gamma\circs \alpha\circs\gamma^{-1})=\md(\gamma)(\mathrm{Inv}(\alpha))$, 
and the surjectivity of the module 
map $\gamma\in \Aut(M)\rightarrow
\md(\gamma)\in \Aut_{\theta}(Z(\tM))$ \cite{Su-Tak-fields}. \hfill$\Box$
\smallskip

Since there exists a genuine action $\beta$ with
$\mathrm{Inv}(\alpha)=\mathrm{Inv}(\beta)$ by \cite[Proposition 1.5.8]{J-act},
\cite[Theorem 5.14]{Su-Tak-act}, \cite[Proposition 22]{KwST},  
the following corollary immediately follows from Theorem
\ref{thm:class}.
(Also see \S \ref{sec:model} about the construction of model actions.)
\begin{cor}\label{cor:2coho}
 The second cohomology vanishing theorem holds for any cocycle crossed action of a discrete
 amenable group on the injective factor of type II$_1$.
\end{cor}
\textbf{Remark.}
So far, it is known that the second cohomology vanishing theorem holds in the
following cases. (See \cite{J-act}, \cite{Ocn-act}, \cite{Su-homoII}.) 
\begin{enumerate}
\itemsep=0pt
\renewcommand{\labelenumi}{(\arabic{enumi})}
 \item arbitrary cocycle crossed actions of free groups on arbitrary von Neumann algebras,
 \item arbitrary cocycle crossed actions of arbitrary locally compact
       groups on properly infinite von Neumann algebras, 
 \item arbitrary cocycle crossed actions of finite groups on type II$_1$ von Neumann algebras,
 \item centrally free cocycle crossed actions of discrete amenable groups on McDuff factors.
\end{enumerate}
In particular, for infinite discrete amenable groups and the injective
factor of type II$_1$, the second cohomology vanishing theorem has been known for
only free cocycle crossed actions. Hence the above corollary removes
the assumption of freeness in this case.  

We close this section by 
giving a cohomological explanation of definition of $\chi(\alpha)$
for a cocycle crossed action $(\alpha,v^\alpha(g,h))$. 
We consider an extended  cocycle crossed
action $(\tal, v^\alpha(g,h))$ of $\tilG$ defined as above. 

We define a multiplication on $\tilU(M)\times \tilG$ by 
$(u,g)(w,h)=(u\tal_g(w)v^\alpha(g,h), gh)$ and denote this group by
$\tilU(M)\rtimes_{v^\alpha}\tilG$. Indeed, the 2-cocycle property of $v^\alpha(g,h)$
assures that this multiplication is associative. 
If an appearing
2-cocycle is clear, then we simply write $\tilU(M)\rtimes\tilG$.
We can easily verify
$(1,g)(1,g^{-1}ng)=(\delv^\alpha(g,n),n)(1,g)$.

Let $\bar{H}:=\{(u^*,n)\mid n\in H \mbox{ and } \Ad u=\tal_n\}$. Then 
$\bar{H}$ is a normal subgroup of $\tilU(M)\rtimes\tilG$. Indeed, 
since we have
$$(u^*,m)(w^*,n)=(u^*\tal_n(w^*)v^\alpha(m,n),mn)
=(w^*u^*v^\alpha(m,n), mn) 
$$ and $\Ad \left(v^\alpha(m,n)^*uw\right)=\Ad v^\alpha(m,n)^*\circs \tal_m\circs\tal_n
=\tal_{mn}$, $\bar{H}$ is a subgroup.
To see $\bar{H}$ is normal, we only have to verify $(w,e)(u^*,n )(w,e)^{-1}$
and $(1,g)(u^*,g^{-1}ng)(1,g)^{-1}$ are in $\bar{H}$.
In fact, 
$$(w,e)(u^*,n)(w,e)^{-1}=(wu^*,n)(w^*,e)=(wu^*\tal_n(w^*),n)=(u^*,n)$$
holds. In particular,  $\tilU(M)$ are in the commutant of $\bar{H}$.

Next we have the following. 
$$(1,g)(u^*,g^{-1}ng)(1,g)^{-1}=(\tal_g(u^*),e)(1,g)(1,g^{-1}ng)(1,g)^{-1}=
(\tal_g(u^*)\delv^\alpha(g,n),n),
$$ 
and 
$$\Ad \left(\delv^\alpha(g,n)^*\tal_g(u)\right)=
\Ad \delv^\alpha(g,n)\circs  \tal_g\circs
\Ad u\,\circs \tal_g^{-1}=
\Ad \delv^\alpha(g,n)\circs  \tal_g\circs
\tal_{g^{-1}ng}\circs \tal_g^{-1}=\tal_n.
$$
Hence $(1,g)(u^*,g^{-1}ng)(1,g)^{-1}\in \bar{H}$.

Although $\{(1,g)\mid g\in G\}$ is not a subgroup of $\tilU(M)\rtimes
\tilG$, $(u^*,n)\rightarrow (1,g)(u^*,n)(1,g)^{-1}$ defines an action of
$\tilG$ on $\bar{H}$ due to $[\tilU(M), \bar{H}]=\{(1,e)\}$.

We have the following $\tilG$-equivariant exact sequence. 
$$0\longrightarrow
 U(Z(\tM))\overset{\iota}{\longrightarrow }\bar{H}
 \overset{\pi}{\longrightarrow } H
\longrightarrow 0.$$
Here $\iota(z)=(z^*,e)\in \bar{H}$, and $\pi(u^*, n)=n$.
 The 
characteristic invariant associated with this exact sequence 
is nothing but $\chi(\alpha)$. 
See
\cite{J-act}, \cite{Su-Tak-act} and \cite{KtST} for cohomological property of
$\chi(\alpha)$. 

\noindent
\textbf{Remark} Let $\tM\rtimes_{\tal, v} \tilG$ be a twisted crossed
product, and $\lambda_g$ the implementing unitary.
The group $\tilU(M)\rtimes G$ is identified with $\{u\lambda_g\mid u\in
\tilU(M),g\in G\}$.

\section{Quasi cocycle crossed actions}\label{sec:cocycle}

Let $\alpha$ and $\beta$ be as in Theorem \ref{thm:class}.
 In this case, we can choose $\tilde{u}_n^\beta\in U(\tM)$ such that  
$$\tbe_g(\tilu_{g^{-1}n g}^\beta)=\lambda(g,n)
\delv^\beta(g,n)\tilu_n^\beta, \,\,
\tilde{u}_m^\beta \tilde{u}_n^\beta=\mu(m,n)v^\beta(m,n)\tilde{u}_{m n}^\beta,\,\,
\theta_t(\tilde{u}_m^\beta)=c_{t}(m)\tilde{u}_m^\beta.$$ 

To explain our idea of proof of Theorem \ref{thm:class}, consider  
the following special
case. Set $Q=G/H$. 
Assume that $v^\alpha(g,h)=v^\beta(g,h)=1$, and 
$G$ is of the form
$G=H\rtimes Q$.  
Then $\alpha|_Q$ and $\beta|_Q$ are centrally free actions of $Q$ on $M$
with $\alpha_p\circs \beta_p^{-1}\in \overline{\Int}(M)$. By
\cite{M-inter-mathscan}, $\Ad v_p\circs\alpha_p=\sigma\circs\beta_p\circs\sigma^{-1}$ for
some $\alpha$-cocycle $v_p$ and $\sigma\in \overline{\Int}(M)$.
Put $v_{h}=\sigma(\tilu^\beta_h)
\tilu_h^{\alpha*}$ for $h\in H$. 
By the choice of $\tilu^\alpha_h$ and $\tilu^\beta_h$, we can  verify
that 
 $v_h\in \tM^\theta=M$, and $v_h$
is a 1-cocycle for $\alpha|_{H}$ with $\Ad v_h \circs\alpha_h=\sigma
\circs \beta_h\circs \sigma^{-1}$.
We then define $v_{hp}=v_h\alpha_h(v_p)$ for $h\in H$, $p\in Q$. We can verify
$v_{p}\alpha_p(v_h)=v_{ph}$  for $p\in Q$ and $h\in H$ as follows. 
\begin{eqnarray*}
 v_p\alpha_p(v_h)&=& v_p\tal_p(\sigma(\tilu_h^\beta)\tilu_h^{\alpha*}) = 
\sigma\circs \tbe_p(\tilu_h^\beta)
v_p \tal_p(\tilu_h^{\alpha*}) \\ 
&=& 
\sigma(\lambda(p, php^{-1})\tilu_{php^{-1}}^\beta)
v_p \lambda(p,php^{-1})^*\tilu_{php^{-1}}^{\alpha*} \\ 
&=& 
\sigma(\tilu_{php^{-1}}^\beta)
v_p \tilu_{php^{-1}}^{\alpha*}  
= 
\sigma(\tilu_{php^{-1}}^\beta)\tilu_{php^{-1}}^{\alpha*}  
\alpha_{php^{-1}}(v_p ) \\
&=&v_{php^{-1}}\alpha_{php^{-1}}(v_p )=v_{ph}.
\end{eqnarray*}
It follows that
$v_{g}$ becomes a 1-cocycle for $\alpha$ and    
$\Ad v_g\circs\alpha_g=\sigma\circs \beta_g\circs \sigma^{-1}$, $g\in G$. 
Thus $\alpha$ and
$\beta$ are strongly cocycle conjugate. 

We would like to extend the above argument to general case. 
Main difficulty is that $G$ is not a
semidirect product of $H$ by $Q$, that is, we can not embed $Q$ into $G$
as a subgroup.
So $\alpha$ does not give an action of $Q$. To treat
such case, we introduce the notion of quasi cocycle crossed cocycle actions of $Q$.
We first classify two quasi cocycle actions of $Q$ by intertwining
argument, 
and combine $H$-parts and $Q$-parts as above.  In what follows, we
mainly use letters $g,h,k$ for elements of $G$, $m,n$ for those of $H$, 
and $p,q,r$ for
those of $Q$.

\begin{lem}\label{lem:innervanish}
 Define $u_m=\tilde{u}_m^\alpha \tilde{u}_m^{\beta*}$. Then 
$u_m\in M$ and we have  $$\alpha_m=\Ad u_m\circs\beta_m,\,\,v^\alpha(m,n)=u_m\beta_m(u_n)
v^\beta(m,n)u_{m n}^*.$$ 
\end{lem}
\textbf{Proof.} By the choice of $\tilde{u}_m^\alpha$ and
 $\tilde{u}_m^\beta$,
 $\theta_t(u_m)=u_m\in \tM^\theta =M$.
It is clear $\alpha_m=\Ad u_m\circs \beta_m$, since $\tal_m=\Ad
u_m \circs \tbe_m$.
The last equation can be verified as follows.  
\begin{eqnarray*}
 u_m\beta_m(u_n)v^\beta(m,n)u_{m n}^*&=&\tilu_m^\alpha\tilu_m^{\beta*}
\tbe_m(\tilu_n^\alpha\tilu_n^{\beta*})v^\beta(m,n)\tilu_{m n}^\beta\tilu_{m
n}^{\alpha*} \\
&=&\tilu_m^\alpha
\tilu_n^\alpha\tilu_n^{\beta*}\tilu_m^{\beta*}
v^\beta(m,n)\tilu_{m n}^\beta\tilu_{m n}^{\alpha*} \\
&=&\mu(m,n)^*\tilu_m^\alpha
\tilu_n^\alpha\tilu_{m n}^{\alpha*} \\
&=&v^\alpha(m,n).
\end{eqnarray*}
\hfill$\Box$

By the above lemma, we can assume $(\alpha_n, v^\alpha(m,n))=(\sigma_n,
v^\sigma(m,n))$ for some fixed cocycle crossed action $\sigma$ of $H$ by a
suitable unitary perturbation. 
In fact, the existence of model actions allows us to further 
assume $v^\sigma(m,n)=1$ for all $m,n\in H$. In what follows, we fix
$\tilu_n\in \tilU(M)$ with $\tilde{\sigma}_n=\Ad \tilu_n$.

\smallskip

\noindent
\textbf{Remark.} Let $\gamma_g$ be a unitary perturbation of
$\alpha$ by $w_g\in U(M)$. If we choose $\tilu^{\gamma}_n$ as $w_n\tilu^\alpha_n$,
then we obtain a same characteristic cocycle.

\smallskip

Fix a section $p\in Q\rightarrow \secp\in G$ with $\tilde{e}=e$, and set
$m(p,q)=\secpq^{-1}\secp\sq\in H$.  So $\secp\sq=\secpq\cdot m(p,q)$ holds.
From $(\secp\sq)\sr=\secp(\sq\sr)$, we obtain
$m(pq,r)\sr^{-1}m(p,q)\sr=m(p,qr)m(q,r)$. We  denote this
element by $m(p,q,r)$.

Set 
$c^\alpha(p,q)=v^\alpha(\secp,\sq)
v^{\alpha}(\secpq,m(p,q))^*$.   
It is easy to see ${\alpha}_{\secp}\circs{\alpha}_{\sq}=\Ad
c^\alpha(p,q)\circs{\alpha}_{\secpq}\circs\sigma_{m(p,q)}$. 

\begin{df}\label{df:q-2-coho}
We say 
$({\alpha}_{\secp}, c^\alpha(p,q))$ a quasi
 cocycle crossed action of $Q$.
\end{df}

The unitary $c^\alpha(p,q)$ behaves like a 2-cocycle as follows.

\begin{lem}\label{lem:K-kercocycle}
For $p,q,r \in Q$, we have
$$
 c^\alpha(p,q)
{\alpha}_{\secpq}\left(\delv^\alpha
(\sr,m(p,q))^* \right) 
c^\alpha(pq,r)
={\alpha}_{\secp}(c^\alpha(q,r))c^\alpha(p,qr).
$$
\end{lem}
\textbf{Proof.} 
In $\tilU(M)\rtimes_{v^\alpha}\tilG$, we compute
$\left((1,\secp)(1,\sq)\right)(1,\sr)=(1,\secp)\left((1,\sq)(1,\sr)\right)$. 
First note the following relations.  
$$(1,\secp)(1,\sq)=(v^{\alpha}(\secp,\sq),\secp\sq)=
(c^\alpha(p,q),\secpq)(1, m(p,q)), $$
$$(1,n)(1,g)=(\delv^\alpha(g,n)^*,e)(1,g)(1,g^{-1}ng).$$

On one hand, we have the following.
\begin{eqnarray*}
\lefteqn{ \left((1,\secp)(1,\sq)\right)
(1,\sr)} \\
&=&
(c^\alpha(p,q),\secpq)(1,m(p,q))(1,\sr)\\
&=&
(c^\alpha(p,q),\secpq)\left(\delv^\alpha(\sr,m(p,q))^*,e\right)
(1,\sr)(1,\sr^{-1}m(p,q)\sr)\\
&=&(c^\alpha(p,q),e)
\left(\alpha_{\secpq}\left(\delv^\alpha(\sr,m(p,q))^*\right),e\right)
(1,\secpq)(1,\sr)(1,\sr^{-1} m(p,q)\sr)\\
&=&
\left(c^\alpha(p,q)\alpha_{\secpq}\left(\delv^\alpha(\sr,m(p,q))^*
\right)c^\alpha(pq,r),\secpqr\right)
(1, m(pq,r))(1,\sr^{-1} m(p,q)\sr) \\
&=&
\left(c^\alpha(p,q)\alpha_{\secpq}\left(\delv^\alpha(\sr,m(p,q))^*
\right)c^\alpha(pq,r),\secpqr\right)
(1, m(p,q,r)).
\end{eqnarray*}

On the other hand, we have the following.
\begin{eqnarray*}
(1,\secp)\left((1,\sq)(1,\sr)\right) &=&
(1,\secp)(c^\alpha(q,r),\sqr)(1,m(q,r)) \\ &=&
(\alpha_{\secp}(c^\alpha(q,r)),e)(1,\secp)(1,\sqr)(1,m(q,r)) \\ 
&=&
(\alpha_{\secp}(c^\alpha(q,r))c^\alpha(p,qr),e)(1,\secpqr)
(1,m(p,qr))(1,m(q,r))\\
&=&
(\alpha_{\secp}(c^\alpha(q,r))c^\alpha(p,qr),\secpqr)
(1,m(p,q,r)).
\end{eqnarray*}
Note $v^\alpha(m,n)=1$ for $m,n\in
H$. 
Thus we get the conclusion. 
\hfill$\Box$

\smallskip
\noindent
\textbf{Remark.} If we do not assume $v^\alpha(m,n)=1$, then we get 
\begin{eqnarray*}
&&c^\alpha(p,q)
{\alpha}_{\secpq}\left(\delv^\alpha
(\sr,m(p,q))^* \right) 
c^\alpha(pq,r)\alpha_{\secpqr}(v^\alpha(m(pq,r),\sr^{-1}m(p,q)\sr)) \\
&&=
{\alpha}_{\secp}(c^\alpha(q,r))c^\alpha(p,qr)\alpha_{\secpqr}(v^\alpha(m(p,qr),m(q,r))).
\end{eqnarray*}

\section{Approximation of quasi cocycle crossed actions}\label{sec:ultra}

In the rest of paper, we use the following notation. For $\alpha\in
\Aut(M)$, $a\in M$ 
and $\psi\in M_*$, 
functionals $\alpha(\psi)$, $a\cdot \psi$ and $\psi \cdot a$ are defined
as follows.
$$\alpha(\psi)=\psi\circs \alpha^{-1},\,\, a\cdot \psi(x)=\psi(xa),\,\,
\psi\cdot a(x)=\psi(ax).$$
Hence $\{\alpha_\nu\}_\nu$ converges to $\alpha$ in the $u$-topology if
$\lim_{\nu}\|\alpha_\nu(\psi)-\alpha(\psi )\|=0$ for all $\psi \in
M_*$. 
It is easy to see a norm bounded sequence $\{a_\nu\}_\nu$ converges to $a$ in the
$\sigma$-strong* topology if and only if 
$$\lim_{\nu}\|(a_\nu-a)\cdot
\psi\|=\lim_{\nu}\|\psi\cdot(a_\nu-a)\|=0$$ 
for all $\psi\in 
M_*$, and equivalently 
$$\lim_{\nu}\|(a_\nu-a)\cdot 
\varphi \|=\lim_{\nu}\|\varphi \cdot(a_\nu-a)\|=0$$ 
for some fixed faithful normal state $\varphi$.
The advantage of use
of these norms instead of usual norms $\|a\|_\varphi^\#$ defining
$\sigma$-strong* topology is the unitary invariance of  $\|\psi\|$, 
i.e., $\|u\cdot \psi\cdot v \|=\|\psi\|$ for $\psi\in M_*$, $u,v\in
U(M)$. 

Throughout this paper, we fix a free ultrafilter $\omega$ over
$\mathbb{N}$. Let $M^\omega$ be an ultraproduct algebra \cite{Ocn-act},
and $M_\omega$ a central sequence algebra \cite{Mc}, \cite{Co-almost}. 
By the Connes-Krieger-Haagerup classification of injective factors
\cite{Co-inj}, \cite{Kri-erg}, \cite{Co-III1}, \cite{Ha-III1}, $M$ is a
McDuff factor. Hence $M_\omega$ is of type II$_1$. 
Let $\tau^\omega :M^\omega\rightarrow M$ be a positive map
$\tau^\omega(X)=\lim_{\nu\rightarrow \omega} x_\nu$, $X=(x_\nu)\in
M^\omega$, where limit is taken in the $\sigma$-weak topology. We extend
$\varphi\in M_*$ to $(M^\omega)_*$ by $\varphi(X)=\varphi\circs
\tau^\omega(X)$. The restriction $\tau_\omega=\tau^\omega|M_\omega$
gives a tracial state on $M_\omega$. 
We denote the $L^1$-norm $\tau_\omega(|X|)$ on $M_\omega$ by $|X|_1$. 
We write $A\Subset B$ if $A$ is a finite subset of $B$.

Let $\alpha$ and $\beta$ be as in the previous section.
By \cite[Theorem 1]{KwST}, $\mathrm{Ker}(\md)=\overline{\Int}(M)$, and
$\Cnt_r(M)=\Cnt(M)$ for an injective factor $M$, where $\Cnt(M)$
is the set  of centrally trivial automorphisms.
Hence 
there exists a sequence of unitaries
$\left\{u_{\sr}^\nu\right\}_\nu\subset  M$, $r\in Q$, such that ${\alpha}_{\sr}=\lim_{\nu}
u_{\sr}^\nu{\beta}_{\sr}$, and 
${\alpha}_{\sr}$ and ${\beta}_{\sr}$ induce
free actions of $Q$ on $M_\omega$. 
However to apply the intertwining argument in our setting, we
also need to approximate 2-cocycles of $\alpha$ and $\beta$. 
This section is devoted to solve
this problem. Our goal in this section is Lemma \ref{lem:appro}.

Let $\{u_{\sr}^\nu\}\subset M$ be as above, and set $U_{\sr}=(u_{\sr}^\nu)\in M^\omega$.
We have $\alpha_{\sr}=\Ad U_{\sr}\circs \beta_{\sr}$ on $M$. 
We want to extend $U_{\sr}$ for all $g\in G$ so that 
$\alpha_g=\Ad U_g\circs \beta_g$ on $M$. It is easy to see two unitaries 
$v^\alpha(n,\sr)\sigma_n(U_{\sr})v^\beta(n,\sr)^{*}$ and 
$v^\alpha(\sr,\sr^{-1}n\sr)^*U_{\sr}v^\beta(\sr,\sr^{-1}n\sr)$ have the
desired property for $\alpha_{n \sr }$ and $\beta_{n \sr }$.
The following lemma says that these two unitaries coincide.



\begin{lem}\label{lem:cocycle0}
For any $W_g=(w_g^\nu)\in M^\omega$ with $\alpha_g=\lim_{\nu\rightarrow \omega}w_g^\nu\circs \beta_g$, 
we have $$v^{\alpha}(n,g)^*\sigma_n(W_{g})v^\beta(n,g)=
v^{\alpha}(g, g^{-1}ng)^*W_{g}v^{\beta}(g, g^{-1}ng)$$ for $g\in G$ and
 $n\in H$.
\end{lem}
\textbf{Proof.}  
We will show
$W_{g}^*\delv^{\alpha}(g, n)\sigma_n(W_{g})=
\delv^{\beta}(g, n)$.
Since the canonical extension is continuous in the $u$-topology, 
we have $\lim\limits_{\nu\rightarrow \omega}\Ad
w_g^\nu\circs\tilde{\beta}_{g}=\tilde{\alpha}_{g}$. Then it follows that
\begin{eqnarray*}
\lim_{\nu\rightarrow\omega }
 w_{g}^{\nu*}\delv^{\alpha}(g, n)
\sigma_n(w_{g}^\nu)
&=& \lim_{\nu\rightarrow\omega }
w_{g}^{\nu*}\delv^{\alpha}(g, n)
\tilde{\sigma}_n(w_{g}^\nu) 
= \lim_{\nu\rightarrow\omega } 
w_{g}^{\nu*}\delv^{\alpha}(g, n)
\tilu_n w_{g}^\nu \tilu_n^* \\
&=& 
\tbe_{g}\circs\tal_{g}^{-1}\left(
\delv^{\alpha}(g, n)
\tilu_n \right)\tilu_n^* 
= \lambda(g,n)^*
\tbe_{g}\left(
\tilu_{g^{-1}ng}\right) \tilu_n^* \\
&=& \delv^\beta(g,n).
\end{eqnarray*}
Here note $\tbe_{g}\circs \tal_{g}^{-1}(\lambda(g,n))=\lambda(g,n)$,
since $\md(\alpha_{g})=\md(\beta_{g})$.
\hfill$\Box$
\medskip

Define $U_{n\sr }=v^\alpha(n,\sr)\sigma_n(U_{\sr})v^\beta(n,\sr)^{*}$, 
$\gamma_g=\Ad U_g\circs {\beta}_g$,
 $V(g,h)=U_g\beta_g(U_h)v^\beta(g,h)U_{g h}^*$. Then $(\gamma_g, V(g,h))$
 is a cocycle crossed action on $M^\omega$ with $\gamma_g|_M=\alpha_g$. 
Note $V(m,n)=1$, $m,n\in H$, since 
$U_n=1$ for $n\in H$. 
Our first task is to show that ${\gamma}_{\secp}$ is a free cocycle crossed action of
$Q$ on $M_\omega$.

\begin{lem}\label{lem:cocycle}
For $g\in G$ and $n\in H $, we have 
$V(n,g)=v^\alpha(n,g)$ and $V(g,n)=v^\alpha(g,n)$.
\end{lem}
\textbf{Proof.} 
From the definition of $U_{n\sr}$, it follows that 
$$V(n,\sr)=\sigma_n(U_{\sr})v^\beta(n,\sr)U_{n\sr}^*
=\sigma_n(U_{\sr})v^\beta(n,\sr)v^{\beta }(n,\sr)^*\sigma_n(U_{\sr}^*)v^\alpha(n,\sr)
=v^\alpha(n,\sr). $$
By Lemma \ref{lem:cocycle0}, we have 
$U_{\sr n}=U_{\sr n \sr^{-1} n}=
v^{\alpha}(\sr,n)^*U_{\sr}v^\beta(\sr,n)$. Hence
$$V(\sr,n)
=U_{\sr}v^\beta(\sr,n)U_{\sr n}^*
=U_{\sr}v^\beta(\sr,n)v^{\beta}(\sr,n)^*U_{\sr}^*v^\alpha(\sr,n)
=v^\alpha(\sr,n)$$
holds.
If $g= m\sr$, then
\begin{eqnarray*}
V(n,g)&=&V(n,m\sr)=\sigma_n(V(m,\sr)^*)V(n,m)V(n m,\sr)\\ &=&
\sigma_n(v^{\alpha}(m,\sr)^*)v^\alpha(n m,\sr) =v^\alpha(n,m\sr)=v^\alpha(n,g)
\end{eqnarray*}
holds. In a similar way as above, we can verify $V(\sr m,n)=
v^\alpha(\sr m, n)$. 
\hfill$\Box$

\begin{lem}
 Put $z(p,q)=U_{\secp}{\beta}_{\secp}(U_{\sq})
c^\beta(p,q)U_{\secpq}^*c^{\alpha}(p,q)^*$, $p,q\in Q$. 
Then $z(p,q)\in M_\omega$.
\end{lem}
\textbf{Proof.}
For two sequences $\{\psi_\nu\}$, $\{\varphi_\nu\}\subset M_*$, we write
$\psi_\nu \sim \varphi_\nu$ if $\lim\limits_{\nu\rightarrow \omega}\|\psi_\nu-\varphi_\nu\|=0$.
Then \begin{eqnarray*}
 {\alpha}_{\secp}\circs {\alpha}_{\sq}(\varphi)
&\sim&
\Ad u_{\secp}^\nu{\beta}_{\secp}(u_{\sq}^\nu) \circs{\beta}_{\secp}\circs{\beta}_{\sq}(\varphi) \\
&=&
\Ad \left(u_{\secp}^\nu{\beta}_{\secp}(u_{\sq}^\nu)c^\beta(p,q) \right)
\circs {\beta}_{\secpq}\circs \sigma_{m(p,q)}(\varphi) \\
&\sim&
\Ad \left(u_{\secp}^\nu{\beta}_{\secp}(u_{\sq}^\nu)c^\beta(p,q) 
u_{\secpq}^{\nu*}\right)\circs 
{\alpha}_{\secpq}\circs \sigma_{m(p,q)}(\varphi) \\
&=&
\Ad \left(u_{\secp}^\nu{\beta}_{\secp}(u_{\sq}^\nu) c^\beta(p,q)
u_{\secpq}^{\nu*}c^\alpha(p,q)^*\right)\circs 
{\alpha}_{\secp}\circs {\alpha}_{\sq}(\varphi) 
\end{eqnarray*}
holds. This implies
$U_{\secp}{\beta}_{\secp}(U_{\sq})c^\beta(p,q)U_{\secpq}^*
c^\alpha(p,q)^*\in M_\omega$. \hfill$\Box$

\begin{lem}\label{lem:cocycleaction}
$(\gamma_{\secp},z(p,q))$ is
a cocycle crossed action of $Q$ on $M_\omega$. 
\end{lem}
\textbf{Proof.}
Set $z'(p,q)=V(\secp,\sq)V(\secpq,m(p,q))^*=
U_{\secp}{\beta}_{\secp}(U_{\sq})c^\beta(p,q)U_{\secpq}^*$. 
Then $(\gamma_{\secp},z'(p,q))$ is a quasi cocycle crossed action of $Q$, and 
$z'(p,q)=z(p,q)c^\alpha(p,q)$.

For $x\in M_\omega$, we have
$$\gamma_{\secp}\circs\gamma_{\sq}(x)=\Ad z'(p,q)\circs
\gamma_{\secpq}\circs
\sigma_{m(p,q)}(x)=\Ad
z(p,q)\circs \gamma_{\secpq}(x),$$ 
since $\Ad c^\alpha(p,q), \sigma_{m(p,q)}\in \Cnt(M)$.

We next show that $z(p,q)$ satisfies the 2-cocycle relation.
By Lemma \ref{lem:K-kercocycle}, we have 
$$
 z'(p,q)\gamma_{\secpq}\left(\delta V(\sr, m(p,q))^*\right)z'(pq,r)
=
\gamma_{\secp}(z'(q,r))z'(p,qr).$$

By Lemma \ref{lem:cocycle} and $\gamma_g|_M=\alpha_g$, the left hand side is 
\begin{eqnarray*}
\lefteqn{ 
z'(p,q)\gamma_{\secpq}\left(\delta V\left(\sr,m(p,q)\right)^*\right)z'(pq,r)}\\ 
&=&
 z(p,q)c^\alpha(p,q){\alpha}_{\secpq}\left(
\delv^\alpha(\sr,m(p,q))^*\right)z(pq,r)c^\alpha(pq,r)\\
&=&
 z(p,q)z(pq,r)c^\alpha(p,q){\alpha}_{\secpq}\left(
\delv^{\alpha}(\sr,m(p,q))^*\right)c^\alpha(pq,r).
\end{eqnarray*}

The right hand side is 
\begin{eqnarray*}
\gamma_{\secp}(z'(q,r))z'(p,qr)&=&
\gamma_{\secp}(z(q,r)c^\alpha(q,r))z(p,qr)c^\alpha(p,qr) \\
&=&
\gamma_{\secp}(z(q,r))z(p,qr){\alpha}_{\secp}(c^\alpha(q,r))c^\alpha(p,qr). 
\end{eqnarray*}
Then Lemma \ref{lem:K-kercocycle} yields
$z(p,q)z(pq,r)=\gamma_{\secp}(z(q,r))z(p,qr)$. 
\hfill $\Box$

\medskip
Next step is to replace $U_p$ so that $z(p,q)=1$. To this end, we
need the second cohomology vanishing theorem. 

\begin{prop}\label{prop:2coho}
Let $Q$ be a discrete amenable group, and 
$(\gamma,u(g,h))$ a semiliftable cocycle crossed action of $Q$ on
 $M_\omega$. Then $u(g,h)$ is a coboundary.
\end{prop}

Proposition \ref{prop:2coho} was first proved by Ocneanu \cite{Ocn-act}, and
later more simplified proof was given in  \cite{M-T-CMP} with
generalization to discrete amenable Kac algebra case. In appendix, 
we present a proof of Proposition \ref{prop:2coho} based on the argument
in \cite{M-T-CMP} for readers' convenience.

\begin{lem}\label{lem:appro}
 Let $\alpha$ and $\beta$ be as above. Then there exists a sequence 
 $\{u_p^\nu\}_{\nu} \subset U(M)$, $p\in Q$, such that \\
$(1)$ $\alpha_{\secp}=\lim_{\nu\rightarrow \omega} 
\Ad u_p^\nu \circs\beta_{\secp}$ in the $u$-topology, \\
$(2)$ $c^\alpha(p,q)=\lim_{\nu\rightarrow \omega}u_p^\nu
\beta_{\secp}(u_q^\nu)c^\beta(p,q)u_{pq}^{\nu*}$ in the  $\sigma$-strong* topology.
\end{lem}

\noindent
\textbf{Proof.} Let $(\gamma_{\secp}, z(p,q))$ be as in Lemma \ref{lem:cocycleaction}.
By Proposition \ref{prop:2coho}, 
there exists $a_p\in U(M_\omega)$ such that
$a_p\gamma_{\secp}(a_q)z(p,q)a_{pq}^*=1
$. Then we have 
\begin{eqnarray*}
1&=&a_p\gamma_{\secp}(a_q)z(p,q)a_{pq}^*=a_pU_{\secp}\beta_{\secp}(a_q)U_{\secp}^*
U_{\secp}{\beta}_{\secp}(U_{\sq})c^\beta(p,q)U_{\secpq}^*c^\alpha(p,q)^*a_{pq}^*\\ &=&
a_pU_{\secp}\beta_{\secp}(a_qU_{\sq})c^\beta(p,q)U_{\secpq}^*a_{pq}^*c^\alpha(p,q)^*.
\end{eqnarray*}
Thus the representing sequence of  
$a_pU_{\secp}$ is a desired one.
\hfill$\Box$

\smallskip

The following lemma is unnecessary in the rest of this paper. However
the similar argument will be appear in the proof of main theorem in
\S\ref{sec:inter}.  

\begin{lem}\label{lem:2cocycle}
Let $U_{\secp}=(u_p^\nu)\in M^\omega$ be as in  Lemma \ref{lem:appro}. 
Then we have $V(g,h)=v^\alpha(g,h)$ for all $g,h\in G$. 
 \end{lem}
\textbf{Proof.} 
By Lemma \ref{lem:appro}, we have 
$$v^\alpha(\secp, \sq)v^\alpha(\secpq, m(p,q))^*=c^\alpha(p,q)=
U_{\secp}{\beta}_{\secp}(U_{\sq})c^\beta(p,q)U_{\secpq}^*=V(\secp, \sq)V(\secpq,m(p,q))^*.
$$
Then 
by Lemma \ref{lem:cocycle}, 
we have
$V(\secp,\sq)=v^\alpha(\secp,\sq)$. 
If $g=m\secp$, then 
\begin{eqnarray*}
 V(g,\sq)&=&V(m\secp,\sq)=V(m,\secp)^*\sigma_m\left(V(\secp,\sq)\right)V(m,\secp\sq) \\ &=&
v^\alpha(m,\secp)^*\sigma_m(v^\alpha(\secp,\sq))v^\alpha(m,\secp\sq)
=v^\alpha(m\secp,\sq)=v^\alpha(g,\sq) 
\end{eqnarray*}
holds. For $h=n\sq$, we have
\begin{eqnarray*}
V(g,h)&=&V(g,n\sq)=\gamma_g(V(n,\sq)^*)V(g,n)V(gn,\sq)\\ &=&
\alpha_g(v^\alpha(n,\sq)^*)v^\alpha(g,n)v^\alpha(gn,\sq)=v^\alpha(g,n\sq)=v^\alpha(g,h) 
\end{eqnarray*}
by Lemma \ref{lem:cocycle} and the above result. 
\hfill$\Box$

\section{Approximate cohomology vanishing}\label{sec:coh}
Let $K$ be a discrete amenable group. Let  $F\Subset K$ and $\delta>0$. 
In this paper, we say $S\Subset K$ is $(F,\delta)$-invariant if 
$$\left|S\cap \bigcap\limits_{g\in F}g^{-1}S\right|>(1-\delta)|S|.$$

We will use the following Rohlin type theorem to show approximate
cohomology vanishing.
\begin{thm}\label{thm:Rohlin}
Let $K$ be a discrete amenable group, and $\gamma_g$ a semiliftable
 strongly free action of $K$ on $M_\omega$. Fix $e\in F\Subset K$,
 $\delta>0$, and  let $S\Subset K$ be an $(F,\delta)$-invariant finite
 set. Then there exists a partition of unity $\{E_s\}_{s\in S}\subset
 M_\omega$ 
such  that
$$\sum_{s\in g^{-1}S\cap S}\left|\gamma_g(E_s)-E_{gs}\right|_1< 4\delta^{\frac{1}{2}},
\,\,\, g\in F,$$
$$\sum_{s\in  S\backslash g^{-1}S}|E_s|_1< 3\delta^{\frac{1}{2}},\,\,\,g\in F,$$
$$\left[\gamma_g(E_s),E_{s'}\right]=0, \,\,\, g\in F,\,\, s,s'\in S.$$
Moreover we can take $\{E_s\}$ in the relative commutant of any
 countable subset of $M_\omega$.
\end{thm}
The proof of Theorem \ref{thm:Rohlin} is essentially same as that of
\cite[Theorem 6.1]{Ocn-act}. (See Appendix for detail of proof.)
Since we do not need a paving structure of a
discrete amenable group, our Rohlin type theorem takes a simpler form
than Ocneanu's one.

The following lemma is the key for our intertwining argument in
the next section.
\begin{lem}\label{nearvanish}
Let $K$ be a discrete amenable group.
Assume $F\Subset K $, $\Psi\Subset M_*$, $\Phi\Subset (M_*)_+$ and
 $\varepsilon>0$ are given.  
Let $S$ be an $(F,\epsilon)$-invariant finite set. 
Let $\gamma$ be a map from $K$ into  
 $\Aut(M)$ such that 
$\gamma_{gh}\equiv \gamma_g\gamma_h$, and
$\gamma_g\not\equiv  \id $ for $g\ne e$, modulo $\Cnt(M)$. 
Assume that a family of unitaries $\{u_g\}_{g\in K}\subset U(M)$ satisfies  
$$\|[\psi, u_s]\|<(3|S|)^{-1}\epsilon,\,\, s\in S,\,\psi\in \Psi,$$
$$\|\varphi\cdot(u_g\gamma_g(u_s)u_{gs}^*-1)\|<5\sqrt{\epsilon}, \,\,
\|(u_g\gamma_g(u_s)u_{g s}^*-1)\cdot \varphi\|<5\sqrt{\epsilon}, \,\, 
\varphi\in \Phi,\, g\in F, \,s\in S.$$
Then there exists
 $w\in U(M)$ such that $$\|[w,\psi]\|<\epsilon,\,\,\,  \psi\in \Psi,$$
$$\|(u_g\gamma_g(w)w^*-1)\cdot \varphi \|<7\sqrt[4]{\epsilon},\,
\|\varphi\cdot(u_g\gamma_g(w)w^*-1) \|<7\sqrt[4]{\epsilon},\,g\in F,\, \varphi\in \Phi.$$
\end{lem}
\textbf{Proof.} By assumption, $\gamma$ induces a free action of $K$ on
$M_\omega$. 
Let $\{E_s\}_{s\in S}\subset M_\omega$ be Rohlin
projections for $\gamma_g$ as in Theorem \ref{thm:Rohlin}. 
Set $W=\sum_{h\in S}u_h E_h\in U(M^\omega)$. 
Note $|\sum_ix_i p_i|=\sum_i|x_i|p_i$ for $x_k\in M$ and a
partition of unity $\{p_i\}\subset M_\omega$,
and $\varphi(ab)=\varphi(a)\tau_\omega(b)$ for $a\in M$ and $b\in M_\omega$.
Then we have
\begin{eqnarray*}
 \varphi \left(\left|u_g\gamma_g(W)W^*-1\right|\right)
&=&
 \varphi \left(\left|\sum_{h,k\in S}(u_g\gamma_g(u_h)u_{k}^*-1)\gamma_g(E_h)E_k\right| \right)\\
&=&\sum_{h,k\in S}
 \varphi \left(\left|u_g\gamma_g(u_h)u_{k}^*-1\right| \gamma_g(E_h)E_k
\right) \\
&=&\sum_{h, k\in S}
 \varphi \left(\left|
u_g\gamma_g(u_h)u_{k}^*-1\right|
 \right)\left|\gamma_g(E_h)E_k\right|_1.
\end{eqnarray*}
We divide $\sum_{h,k\in S}$ as $$\sum_{h,k\in S}=
\sum_{\stackrel{\scriptstyle h\in S\cap 
g^{-1}S}{k\in S}}
+\sum_{\stackrel{\scriptstyle h\in S\backslash g^{-1}S}{k\in S}}=
\sum_{\stackrel{\scriptstyle h\in S\cap g^{-1}S}{g h=k\in S}}+
\sum_{\stackrel{\scriptstyle h\in S\cap 
g^{-1}S}{g h\ne k \in S}}+
\sum_{\stackrel{\scriptstyle h\in S\backslash g^{-1}S}{k\in S}}.
$$
Since $\varphi(|x|)\leq \|x\cdot\varphi\|$, the first term is less than 
$5\sqrt{\epsilon}$ by the assumption.
The second term is estimated as follows.
\begin{eqnarray*}
\sum_{\stackrel{\scriptstyle h\in S\cap g^{-1}S}{k\in S,k\ne g h}}
 \varphi \left(\left|
(u_g\gamma_g(u_h)u_{k}^*-1)\right| \right)|\gamma_g(E_h)E_k|_1 
&\leq &2\sum_{h\in S\cap g^{-1}S}|(1-E_{g h})\gamma_g(E_h)|_1 \\
&=&2\sum_{h\in S\cap g^{-1}S}|(1-E_{g h})(\gamma_g(E_h)-E_{g h})|_1 \\
&< & 8\sqrt{\epsilon}.
\end{eqnarray*}
The third term is estimated as follows
$$
 \sum_{\stackrel{\scriptstyle h\in S\backslash g^{-1}S}{k\in S}} \varphi \left(\left|
u_g\gamma_g(u_h)u_{k}^*-1\right|\right)\left|\gamma_g(E_h)E_k\right|_1\leq 
 \sum_{h\in S\backslash g^{-1}S} 2|\gamma_g(E_h)| 
<6\sqrt{\epsilon}.
$$
Hence we obtain
$\varphi\left(\left|u_g\gamma_g(W)W^*-1\right|\right)<19\sqrt{\epsilon}$,
$g\in F$.
In the same way, we have
 $ \varphi\left(\left|(u_g\gamma_g(W)W^*-1)^*\right|\right)<19\sqrt{\epsilon}$.

Let $W=(w_\nu)_\nu$, $E_s=(e_{s,\nu})_\nu$
 be representing sequences consisting of unitaries, and projections respectively.
Set $a_\nu=\sum_{s\in S}u_s e_{s,\nu}$. Note that $a_\nu$ is not a
 unitary in general. Since $(a_\nu)=(w_\nu)$ in $M^\omega$, 
$a_\nu-w_\nu$ converges to $0$
 in the $\sigma$-strong* topology as $\nu\rightarrow \omega$. 

Fix sufficiently large $\nu$
such that $$\|[\psi, w_\nu-a_\nu]\|<\frac{\epsilon}{3},\,\,\|[\psi,
 e_{s,\nu}]\|<\frac{\epsilon}{3|S|},\,\, \psi \in \Psi, s\in S,$$
$$\varphi(\left|u_g\gamma_g(w_\nu)w_\nu^*-1\right|)<19\sqrt{\epsilon}, \,\,
\varphi(\left|(u_g\gamma_g(w_\nu)w^*_\nu-1)^*\right|)<19\sqrt{\epsilon},\,\,g
\in F,\,\, \varphi\in \Phi.$$
It follows that
\begin{eqnarray*}
\left\|[\psi,w_\nu]\right\|&\leq &\|[\psi, w_\nu-a_\nu]\| +
\left\|[\psi,a_\nu]\right\| < \frac{\varepsilon}{3}+
\left\|\left[\psi,\sum_{k\in S}u_ke_{k,\nu}\right]\right\| \\ &\leq&  \frac{\varepsilon}{3}+
\sum_{k\in S}\left(\left\|[\psi,u_k]e_{k,\nu}\right\|+\left\|u_k[\psi,e_{k, \nu}]\right\|\right)   
<\varepsilon.
\end{eqnarray*}
Since $\|x\cdot \varphi \|\leq \sqrt{\|x\|\varphi(|x|)} $ and $\|\varphi\cdot x\|\leq
\sqrt{\|x\|\varphi(|x^*|)}$, $w=w_\nu$ is a desired unitary. \hfill$\Box$

\section{Intertwining argument}\label{sec:inter}
In this section,  we present a proof of the main theorem, Theorem \ref{thm:class}. 
We can assume $\beta$ is
an action, i.e., $v^\beta(g,h)=1$, $g,h\in G$ due to the existence of
model actions. 
Recall that we assumed that 
$\sigma_n=\alpha_n=\beta_n$ is an action of $H$, 
and fixed a unitary $\tilu_m\in
\tM$ with $\tilde{\sigma}_m=\Ad \tilu_m$ in \S \ref{sec:cocycle}.

At first, we will classify two quasi cocycle crossed
actions $(\alpha_{\secp}, c^\alpha(p,q))$ and $(\beta_{\secp}, 1)$ of $Q$.
For simplicity, we write $\alpha_{\secp}$ as $\alpha_p$ until the end
of proof of Theorem \ref{thm:quasiclass}.

\begin{thm}\label{thm:quasiclass}
 Let $(\alpha_p,c^\alpha(p,q))$ and $(\beta_p,1)$ as above. Then there
 exist $\bar{\theta}_0$, $\bar{\theta}_1\in \overline{\Int}(M)$,
 $\hat{u}_p^0$, $\hat{u}_p^1\in U(M)$ such that 
$$\Ad \hat{u}_p^0\circs\bar{\theta}_0\circs \alpha_p\circs
 \bar{\theta}_0^{-1}=
\Ad \hat{u}_p^1\circs\bar{\theta}_1\circs \beta_p\circs
 \bar{\theta}_1^{-1},
$$ and 
\begin{eqnarray*}
 \lefteqn{
\hat{u}_p^0\bar{\theta}_0\circs\alpha_p\circs \bar{\theta}_0^{-1}(\hat{u}_q^0)
\bar{\theta}_0(c^{\alpha}(p,q))
\bar{\theta}_0\circs\alpha_{pq}\circs\bar{\theta}_0^{-1}(
\bar{\theta}_0(\tilu_{m(p,q)})\tilu_{m(p,q)}^*) \hat{u}_{pq}^{0*}} \\ &&=
\hat{u}_p^1\bar{\theta}_1\circs \beta_p\circs\bar{\theta}_1^{-1}(\hat{u}_q^1)
\bar{\theta}_1\circs\beta_{pq}\circs\bar{\theta}_1^{-1}(
\bar{\theta}_1(\tilu_{m(p,q)})\tilu_{m(p,q)}^*) \hat{u}_{pq}^{1*}.
\end{eqnarray*}
\end{thm}
We remark that  $\bar{\theta}_i(\tilu_{m(p,q)})\tilu_{m(p,q)}^*$, $i=0,1$, are
indeed in $M$.

\noindent
\textbf{Proof.} 
Put $\epsilon_n=4^{-n}$, $n\in \mathbb{N}$. (Until the end of proof, we
use the letter $n$ to denote elements in $\mathbb{N}$.)  
Fix $F_n \Subset Q$    
and an $(F_n, \epsilon_n)$-invariant set $S_n\Subset Q$ such that 
$e\in F_1$, $F_n\subset F_{n+1}$, 
$\bigcup_n F_n=Q$, $F_n\subset S_n$ and $F_nS_n\subset F_{n+1}$.

Fix a faithful
normal state $\varphi_0$. Let $\{\Psi_n\}_n $ be an increasing sequence
of finite sets of $M_*$  such that 
$\bigcup_{n}\Psi_n$ is total in $M_*$. 

Set $\gamma^{(0)}=\alpha$, $c^0(p,q)=c^\alpha(p,q)$, 
$\gamma^{(-1)}=\beta$, $c^{-1}(p,q)=1$.  
We will construct a family of quasi cocycle crossed
actions $(\gamma^{(n)}_p, c^n(p,q))$,
unitaries $\cheu_p^n$,
$\baru_p^n$, 
$\checku_p^n$, $w_n$, 
 automorphisms $\theta_n$,
and finite sets $\Phi_n\Subset (M_*)_+$, $\Psi_n', \Phi_n'\Subset M_*$ 
satisfying the following conditions.

\begin{enumerate}
\itemsep=0pt
\renewcommand{\labelenumi}{($n$.\arabic{enumi})}
 \item 
$\displaystyle{\Phi_{n}=\{\Ad \checku_p^{n-1}(\varphi_0)\}_{p\in F_{n}}, \,
\Phi_n'=\bigcup_{\stackrel{\scriptstyle p,q\in F_{n}}
{r\in S_n,\varphi\in \Phi_n}}
\left(\gamma_p^{(n-1)}\right)^{-1}\left(\{\varphi\cdot c^{n-1}(q,r),
c^{n-1}(q,r)\cdot \varphi\}\right), }$ 
 \item 
$\displaystyle{\Psi_{n}'=\Psi_{n}\cup \theta_{n-1}(\Psi_{n})\cup \bigcup_{p\in
F_n}\{\checku_p^{n-1}\cdot\varphi_0, \varphi_0\cdot\checku_p^{n-1}\}\cup \Phi_n', }$ 
 \item 
$\baru_p^n=\cheu_p^n\gamma_p^{(n-2)}(w_n)w_n^*,\,\,
\checku_p^n=\baru_p^n w_n \checku_p^{n-2}w_n^*,\,\,
\theta_n=\Ad w_n \circs \theta_{n-2},\,\, n\geq 3,$
 \item $\gamma_p^{(n)}=\Ad \cheu_p^n\circs \gamma_p^{(n-2)}=
\Ad \baru_p^n\circs \Ad w_n\circs\gamma^{(n-2)}_p\circs \Ad w_n^*, $
 \item $
 c^n(p,q)=\cheu_p^n\gamma_p^{(n-2)}(\cheu_q^n) c^{n-2}(p,q)
\cheu_{pq}^{n*},$
 \item $\displaystyle{
\|\gamma^{(n)}_p(\psi)-\gamma^{(n-1)}_p(\psi)\|<\frac{\epsilon_n}{6|S_n|},\,\,
p\in F_{n+1},}$ \\ 
$\displaystyle{\psi\in \bigcup\limits_{q\in F_{n+1}}
  \left(\gamma_q^{(n-2)}\right)^{-1}(\Psi_{n-1}')
\cup\left(\gamma^{(n-1)}_q\right)^{-1}(\Psi_{n}'), }$
\item 
$\renewcommand{\arraystretch}{1.6}
\left\{\begin{array}{ll}
\displaystyle{\|\varphi\cdot \left(c^n(p,q)-c^{n-1}(p,q)\right)\|<\frac{\epsilon_n}{2}}, & \\
\displaystyle{\|\left(c^n(p,q)-c^{n-1}(p,q)\right)\cdot \varphi\|<\frac{\epsilon_n}{2}},  & 
\varphi\in \Phi_{n}\cup \Phi_{n-1},\, p\in F_{n},\, q\in S_n, 
 \end{array}\right.$
 \item  
$\left\{\begin{array}{ll}
\|\varphi \cdot (\baru^{n}_p-1)\|<7\sqrt[4]{\epsilon_{n-1}}, & \\
\|(\baru^{n}_p-1)\cdot \varphi\|<7\sqrt[4]{\epsilon_{n-1}}, &\varphi\in
\Phi_{n-1},\, p\in F_{n-1}, \,n\geq 2,
\end{array}\right.$
\item $\|[w_n, \psi]\|<\epsilon_{n-1},\, \psi\in \Psi_{n-1}',\, n \geq 2.$
\end{enumerate}

\noindent
\textbf{Step 1.} Define $\Phi_1$, $\Phi_1'$ and $\Psi_1'$ as in $(1.1)$ and
$(1.2)$. (Here we set $\checku_p^{0}=1$, $\theta_0=\id$ and 
$\Psi_0=\Phi_0=\emptyset$.) 
By Lemma \ref{lem:appro}, there exists a unitary 
$\cheu^1_p$ such that 
\begin{eqnarray*}
(1.a)&&\!\!\!\!\!\left\|\Ad
\cheu_p^1\circs \gamma_p^{(-1)}(\psi)-\gamma_p^{(0)}(\psi)\right\|<\frac{\epsilon_1}{6|S_1|},
\,\,p\in F_2,\,
\psi\in \bigcup_{r\in F_2} \left(\gamma_r^{(0)}\right)^{-1}(\Psi'_1),\\
(1.b)&&\!\!\!\!\!
\left\|\varphi\cdot
 \left(\cheu_p^{1}\gamma^{(-1)}_p(\cheu_q^1)\cheu_{pq}^{1*}-
c^0(p,q) \right)\right\|<\frac{\epsilon_1}{2}, \,\,
\left\|\left(\cheu_p^{1}\gamma^{(-1)}_p(\cheu_q^1)\cheu_{pq}^{1*}-
c^0(p,q) \right)\cdot \varphi\right\|<\frac{\epsilon_1}{2}, \\
\,
&&p\in F_1,\, q\in S_1,\,\varphi\in \Phi_1. 
\end{eqnarray*}

Set $$w_1=1,\, u_p^1=\baru_p^1=\cheu_p^1,\, 
\theta_1=\Ad w_1,\, \gamma_p^{(1)}=\Ad \cheu_p^1\circs \gamma_p^{(-1)},\,
c^1(p,q)=\cheu_p^{1}\gamma^{(-1)}_p(\cheu_q)\cheu_{pq}^{1*}.$$
Then we obtain $(1.3)$, $(1.4)$ and $(1.5)$. 
The conditions $(1.6)$ and (1.7) follow from $(1.a)$ and $(1.b)$. 
Hence the first step is complete.
\medskip

Suppose we have constructed up to the $(n-1)$-st step. 

\noindent
\textbf{Step \boldmath${n}$.} Define $\Phi_n$, $\Phi_n'$ and $\Psi'_n$
as in $(n.1)$ and $(n.2)$.
By Lemma \ref{lem:appro}, there exists a
unitary $\cheu^n_p$ such that 
\begin{enumerate}
\itemsep=0pt
\renewcommand{\labelenumi}{($n.\alph{enumi}$)}
 \item 
$\displaystyle{\left\|\Ad \cheu_p^n\circs 
\gamma_p^{(n-2)}(\psi)-\gamma_p^{(n-1)}(\psi)\right\|<
\frac{\epsilon_n}{6|S_n|},}$ \\
 $\displaystyle{ p\in F_{n+1},\, 
\psi\in \bigcup_{q\in F_{n+1}} (\gamma_q^{(n-2)})^{-1}(\Psi_{n-1}')
\cup (\gamma_q^{(n-1)})^{-1}(\Psi_n'), }$
 \item 
$\renewcommand{\arraystretch}{1.6}
\begin{array}[t]{l}
\displaystyle{\left\|\varphi\cdot \left(\cheu_p^{n}\gamma^{(n-2)}_p(\cheu_q^n)
c^{n-2}(p,q)\cheu_{pq}^{n*}-
c^{n-1}(p,q) \right)\right\|<\frac{\epsilon_n}{2}},\\
\displaystyle{\left\|\left(\cheu_p^{n}\gamma^{(n-2)}_p
(\cheu_q^n) c^{n-2}(p,q)
\cheu_{pq}^{n*}-
c^{n-1}(p,q) \right)\cdot \varphi\right\|<\frac{\epsilon_n}{2}, } \\
p\in F_n,\, q\in S_n,\,\varphi\in \Phi_n \cup \Phi_{n-1}.
 \end{array}
$ 
\end{enumerate}

By $(n.a)$ and $(n-1.6)$, we get 
$$\left\|\Ad
\cheu_p^n\circs\gamma_p^{(n-2)}(\psi)-\gamma_p^{(n-2)}(\psi)\right\|<
\frac{\epsilon_{n-1}}{3|S_{n-1}|},\,\,\, \psi\in \bigcup_{q\in F_n}
(\gamma_q^{(n-2)})^{-1}(\Psi_{n-1}'),\,\,p\in F_n.$$ 
Hence 
$$\left\|\left[\cheu_p^n,\psi\right]\right\|<
\frac{\epsilon_{n-1}}{3|S_{n-1}|},\,\, \psi\in \Psi_{n-1}', \,\,p\in F_n.$$
By $(n.b)$ and $(n-1.7)$, we have 
$$\left\|\varphi\cdot \left(\cheu_p^{n}\gamma^{(n-2)}_p(\cheu_q^n)
c^{n-2}(p,q)\cheu_{pq}^{n*}-
c^{n-2}(p,q) \right)\right\|<\epsilon_{n-1},$$
$$\left\|\left(\cheu_p^{n}\gamma^{(n-2)}_p(\cheu_q^n)c^{n-2}(p,q)
\cheu_{pq}^{n*}-
c^{n-2}(p,q) \right)\cdot \varphi\right\|<\epsilon_{n-1} 
$$
for $p\in F_{n-1}$, $q\in S_{n-1}$ and $\varphi\in \Phi_{n-1}$.

Since $\Phi_{n-1}\cdot c^{n-2}(p,q)\subset \Psi_{n-1}'$, $p\in F_{n-1}$, $q\in S_{n-1}$, 
 and $F_{n-1}S_{n-1}\subset F_n$, we have 
\begin{eqnarray*}
\lefteqn{ \left\|\varphi\cdot (\cheu_{p}^n\gamma_{p}^{(n-2)}
(\cheu_q^n)\cheu_{pq}^{n*}-1 )\right\|} \\
&=&
\left\|\varphi\cdot \cheu_{p}^n\gamma_{p}^{(n-1)}(\cheu_q^n)-
\varphi\cdot \cheu_{pq}^n\right\| \\
&\leq &
\|[\cheu_{pq}^n,\varphi]\|+
\left\|\varphi\cdot
 \cheu_{p}^n\gamma_{p}^{(n-2)}(\cheu_q^n)-\cheu_{pq}^n\cdot 
\varphi\right\| \\
&=&\epsilon_{n-1}+
\left\|\varphi\cdot \cheu_{p}^n\gamma_{p}^{(n-2)}(\cheu_q^n)
c^{n-2}(p,q)-
\cheu_{pq}^n\cdot  \varphi\cdot c^{n-2}(p,q)\right\| \\
&\leq &\epsilon_{n-1}+
\|[\cheu_{pq}^n,  \varphi\cdot c^{n-2}(p,q)]\|+
\left\|\varphi\cdot \cheu_{p}^n\gamma_{p}^{(n-2)}(\cheu_q^n)c^{n-2}(p,q)-
\varphi\cdot c^{n-2}(p,q)\cheu_{pq}^n\right\| \\
&\leq &2\epsilon_{n-1}+
\left\|\varphi\cdot
 \cheu_{p}^n\gamma_{p}^{(n-2)}(\cheu_q^n)c^{n-2}(p,q)\cheu_{pq}^{n*}
-\varphi\cdot c^{n-2}(p,q)\right\| \\
&<&3\epsilon_{n-1}
\end{eqnarray*}
for $\varphi\in \Phi_{n-1}$, $p\in F_{n-1}$ and $q\in S_{n-1}$.

Since $\left(\gamma_p^{(n-2)}\right)^{-1}(c^{n-2}(p',q)\cdot \varphi)\in
\Psi_{n-1}'$ for 
$p,p'\in F_{n-1}$, $q\in S_{n-1}$ and $\varphi\in \Phi_{n-1}$, 
$$\left\|[\gamma_p^{(n-2)}(\cheu_q),c^{n-2}(p',q')\cdot\varphi]\right\|= 
\left\|[\cheu_q, \left(\gamma_p^{(n-2)}\right)^{-1}\left(
 c^{n-2}(p',q')\cdot \varphi\right)]\right\|<\epsilon_{n-1}
$$ for $p,p'\in F_{n-1}$, $q,q'\in S_{n-1}$ and $\varphi\in \Phi_{n-1}$. 
Thus
\begin{eqnarray*}
\lefteqn{ \left\|(\cheu_{p}^n\gamma_{p}^{(n-2)}(\cheu_q^n)
\cheu_{pq}^{n*}-1 )\cdot \varphi
\right\|} \\ 
&=&
 \left\|\gamma_{p}^{(n-2)}(\cheu_q^{n*})\cheu_{p}^{n*}\cdot \varphi -
\cheu_{pq}^{n*}\cdot \varphi\right\| \\
&\leq &
 \left\|[\gamma_{p}^{(n-2)}(\cheu_q^{n*})\cheu_{p}^{n*},
  \varphi ]\right\|+ 
\left\|
\varphi\cdot \gamma_{p}^{(n-2)}(\cheu_q^{n*})\cheu_{p}^{n*}-
\cheu_{pq}^{n*}\cdot \varphi\right\| \\
&\leq &
 \left\|\gamma_{p}^{(n-2)}(\cheu_q^{n*})[\cheu_{p}^{n*},
  \varphi ]\right\|+ 
 \left\|[\gamma_{p}^{(n-2)}(\cheu_q^{n*}),
  \varphi ]\cheu_{p}^{n*}\right\|+ 
\left\|
\varphi\cdot \gamma_{p}^{(n-2)}(\cheu_q^{n*})\cheu_{p}^{n*}-
\cheu_{pq}^{n*}\cdot \varphi\right\| \\
&\leq& 2\epsilon_{n-1}+
 \left\|\varphi\cdot\gamma_{p}^{(n-2)}(\cheu_q^{n*})\cheu_{p}^{n*} 
-\cheu_{pq}^{n*}\cdot \varphi\right\| \\
&=& 2\epsilon_{n-1}+
 \left\|c^{n-2}(p,q)\cdot\varphi \cdot\gamma_{p}^{(n-2)}(\cheu_q^{n*})\cheu_{p}^{n*} 
-c^{n-2}(p,q)\cheu_{pq}^{n*}\cdot \varphi\right\| \\
&\leq & 2\epsilon_{n-1}+
\left\|\left[
\gamma_{p}^{(n-2)}(\cheu_q^{n*})\cheu_{p}^{n*}, c^{n-2}(p,q)\cdot\varphi
\right]
\right\|\\
&& +\left\|\gamma_{p}^{(n-2)}(\cheu_q^{n*})\cheu_{p}^{n*} c^{n-2}(p,q)\cdot\varphi
-c^{n-2}(p,q)\cheu_{pq}^{n*}\cdot \varphi\right\| \\
&\leq & 2\epsilon_{n-1}+
\left\|\gamma_{p}^{(n-2)}(\cheu_q^{n*})
\left[
\cheu_{p}^{n*}, c^{n-2}(p,q)\cdot\varphi
\right]
\right\|+
\left\|\left[
\gamma_{p}^{(n-2)}(\cheu_q^{n*}), 
c^{n-2}(p,q)\cdot\varphi
\right]\cheu_{p}^{n*}
\right\|\\
&& +\left\|\gamma_{p}^{(n-2)}(\cheu_q^{n*})\cheu_{p}^{n*} c^{n-2}(p,q)\cdot\varphi
-c^{n-2}(p,q)\cheu_{pq}^{n*}\cdot \varphi\right\| \\
&\leq & 4\epsilon_{n-1}+
 \left\|\gamma_{p}^{(n-2)}(\cheu_q^{n*})\cheu_{p}^{n*} c^{n-2}(p,q)\cdot\varphi
-c^{n-2}(p,q)\cheu_{pq}^{n*}\cdot \varphi\right\| \\
&= & 4\epsilon_{n-1}+
 \left\|\left(\cheu_p^n\gamma_{p}^{(n-2)}(\cheu_q^{n})c^{n-2}(p,q)\cheu_{pq}^{n*}
-c^{n-2}(p,q)\right)\cdot\varphi\right\| \\
&<&5\epsilon_{n-1}
\end{eqnarray*}
holds for $p\in F_{n-1}$, $q\in S_{n-1}$ and $\varphi\in \Phi_{n-1}$.

By Lemma \ref{nearvanish}, there exists a unitary $w_n$ such that
$$\|(\cheu_{p}^n\gamma_{p}^{(n-2)}(w_n)w_n^*-1)\cdot \varphi
\|<7\sqrt[4]{\epsilon_{n-1}}, \,\,
\|\varphi \cdot (\cheu_{p}^n\gamma_{p}^{(n-2)}(w_n)w_n^*-1)
\|<7\sqrt[4]{\epsilon_{n-1}}$$ and 
$\|[w_n,\psi ]\|<\epsilon_{n-1}$ for $ \varphi \in \Phi_{n-1}$, 
$p\in F_{n-1}$ and  $\psi \in \Psi_{n-1}'$. Put
$\baru_p^n=\cheu_p^n\gamma_p^{(n-2)}(w_n)w_n^*$. We then obtain
$(n.8)$ and $(n.9)$.

Set 
$$b_p^n=\baru_p^n w_nb_p^{n-2}w_n^*,\,\,\theta_n=\Ad
w_n\circs\theta_{n-2},\,\,
c^n(p,q)=\cheu_p^{n}\gamma^{(n-2)}_p(\cheu_q^n)c^{n-2}(p,q)
\cheu_{pq}^{n*},$$
$$\gamma_p^{(n)}=\Ad \cheu_p^n\circs\gamma_p^{(n-2)}=\Ad \baru_p^n\circs 
\Ad w_n\circs \gamma_p^{(n-2)}\circs \Ad w_n^*.$$  

Then we obtain $(n.3)$, $(n.4)$ and $(n.5)$. From $(n.a)$ and
$(n.b)$, $(n.6)$ and $(n.7)$ follow. Thus the $n$-th step is complete,
and we finished the induction. 
\bigskip

We show $\{\theta_{2n}\}$ and $\{\theta_{2n+1}\}$ converge to
some automorphisms. Fix $n_0\in\mathbb{N}$. 
Take $\psi \in \Psi_{n_0}$.
If $n\geq n_0+1$, then $\psi,\theta_{n}(\psi) \in \Psi_{n+1}'$. By
$(n+2.9)$, we have 
$$\left\|\theta_{n}(\psi)-\theta_{n+2}(\psi)\right\|
=\|[w_{n+2},\theta_{n}(\psi)]\|<\epsilon_{n+1},\,\,
\|\theta_{n}^{-1}(\psi)-\theta_{n+2}^{-1}(\psi)\|=\|[w_{n+2},\psi]\|<
\epsilon_{n+1}.$$ Since $\bigcup_n\Psi_n$ is total in $M_*$,
$\{\theta_{2n}^{\pm 1}(\psi)\}$ and
$\{\theta_{2n+1}^{\pm 1}(\psi)\}$ are Cauchy sequences for all $\psi\in M_*$. Hence limits
$\bar{\theta}_0=\lim_{n}\theta_{2n} $ and  
$\bar{\theta}_1=\lim_{n}\theta_{2n+1} $ exist in the $u$-topology.

Next we show the existence of $\lim_{n}\checku_p^{2n}$ and
$\lim_{n}\checku_p^{2n+1}$. 
Fix $p\in F_{n_0+1}$.
Note $\varphi_0\in \Phi_m\subset \Psi_m'$ for all $m$, and $\varphi_0\cdot
\checku_p^n,  \checku_p^n\cdot\varphi_0\in \Psi_{n+1}' $, $\Ad
\checku_p^n(\varphi_0)\in \Phi_{n+1}$ for $n\geq n_0$.
By $(n+2.8)$ and 
$(n+2.9)$, we have
\begin{eqnarray*}
\left\|\varphi_0\cdot(\checku_p^{n+2}-\checku_p^n) \right\|
&=&
 \|\varphi_0\cdot(\baru_p^{n+2}w_{n+2}\checku_p^{n}w_{n+2}^*
-\checku_p^n) \| \\
 &=&
 \|\varphi_0\cdot(\baru_p^{n+2}
-\checku_p^n w_{n+2}\checku_p^{n*}w_{n+2}^*) \| \\
 &\leq&
 \|\varphi_0\cdot(\baru_p^{n+2}-1)\|+\|
\varphi_0\cdot(1-\checku_p^n w_{n+2}\checku_p^{n*}w_{n+2}^*) \| \\
 &\leq&
 7 \sqrt[4]{\epsilon_{n+1}}+
\|\varphi_0\cdot(w_{n+2}\checku_p^{n}-
\checku_p^n w_{n+2} ) \| \\
 &\leq&
 7\sqrt[4]{\epsilon_{n+1}}+
\|\varphi_0\cdot w_{n+2}\checku_p^{n}- w_{n+2} \cdot\varphi_0\cdot\checku_p^n 
 \| +\|[w_{n+2},\varphi_0\cdot \checku_p^{n}]\|
\\
&\leq & 7\sqrt[4]{\epsilon_{n+1}}+2\epsilon_{n+1}.
\end{eqnarray*}

We next estimate $ \|(\checku_p^{n+2}-\checku_p^n)\cdot \varphi_0\|$. 
By $(n+2.8)$, we have
\begin{eqnarray*}
 \|(\checku_p^{n+2}-\checku_p^n)\cdot \varphi_0 \|
&=&
 \|(\baru_p^{n+2}w_{n+2}\checku_p^n w_n^*
-\checku_p^n)\cdot \varphi_0 \| \\
&\leq &
 \|\baru_p^{n+2}(w_{n+2}\checku_p^{n}w_{n+2}^*-\checku_p^n)\cdot \varphi_0\|+
\|(\baru_p^{n+2}-1)\checku_p^{n}\cdot \varphi_0 \| \\
&< &
 \|(w_{n+2}\checku_p^{n}w_{n+2}^*-\checku_p^n)\cdot \varphi_0\|+7\sqrt[4]{\epsilon_{n+1}}.
\end{eqnarray*}
By $(n+2.9)$, the first term is estimated as follows.
\begin{eqnarray*}
 \|(w_{n+2}\checku_p^{n}w_{n+2}^*-\checku_p^{n})\cdot \varphi_0 \|&= &
 \|(\checku_p^{n}w_{n+2}^*-w_{n+2}^*\checku_p^{n})\cdot \varphi_0 \| \\ 
&\leq &
 \|\checku_p^{n}w_{n+2}^*\cdot \varphi_0-
\checku_p^{n}\cdot \varphi_0\cdot w_{n+2}^* \| +
\|[w_{n+2}^*, \checku_p^{n}\cdot \varphi_0]\|
\\ 
&<& 2\epsilon_{n+1}.
\end{eqnarray*}
 
Hence we obtain 
$ \|(\checku_p^{n+2}-\checku_p^n)\cdot \varphi_0 \|< 
7\sqrt[4]{\epsilon_{n+1}}+2\epsilon_{n+1}. $ 
The above  estimation yields that $\hat{u}_p^{0}=\lim_{n}\checku_{p}^{2n}$
and 
$\hat{u}_p^{1}=\lim_{n}\checku_{p}^{2n+1}$ exist in the
$\sigma$-strong* topology. 

We have $\gamma^{(2n)}_p=\Ad \checku^{2n}_p\circs
 \theta_{2n}\circs \alpha_p \circs
 \theta_{2n}^{-1}$ and
$\gamma^{(2n+1)}_p=\Ad \checku^{2n+1}_p\circs
 \theta_{2n+1}\circs \beta_p \circs
 \theta_{2n+1}^{-1}$ by construction.
Letting $n\rightarrow \infty$, we obtain 
$
\Ad \hat{u}_p^{0}\circs\bar{\theta}_0\circs
\alpha_p\circs \bar{\theta}_0^{-1}=
\Ad \hat{u}_p^{1}\circs\bar{\theta}_1\circs
\beta_p\circs \bar{\theta}_1^{-1}$ by $(n.6)$. 

We will show the convergence of $c^{2n}(p,q)$ and $c^{2n+1}(p,q)$. 
Put $\bar{w}_{2n}=w_{2n} w_{2n-2}\cdots w_2$. Of course $\theta_{2n}=\Ad \bar{w}_{2n}$. 
We can easily verify 
$\cheu_p^{2n}\cheu_p^{2n-2}\cdots \cheu_p^{2}=
b_p^{2n}\bar{w}_{2n}\alpha_p(\bar{w}_{2n}^*)$. Then it follows that
\begin{eqnarray*}
c^{2n}(p,q)
&=&\cheu_p^{2n}\gamma^{(2n-2)}_p(\cheu_q^{2n})
c^{2n-2}(p,q)\cheu_{pq}^{2n*}\\
&=&\cheu_p^{2n}\gamma^{(2n-2)}_p(\cheu_q^{2n})
\cheu_p^{2n-2}\gamma^{(2n-4)}_p(\cheu_q^{2n-2})
c^{2n-4}(p,q)\cheu_{pq}^{2n-2*}
\cheu_{pq}^{2n*}\\
&=&
\cheu_p^{2n}
\cheu_p^{2n-2}\gamma^{(2n-4)}_p(\cheu_q^{2n}\cheu_q^{2n-2})
c^{2n-4}(p,q)\cheu_{pq}^{2n-2*}
\cheu_{pq}^{2n*}.
\end{eqnarray*}
We repeat the above computation and obtain the following.
\begin{eqnarray*}
c^{2n}(p,q)
&=&
\cheu_p^{2n}\cdots 
\cheu_p^{2}\alpha_p(\cheu_q^{2n}\cdots\cheu_q^{2})
c^\alpha(p,q)\cheu_{pq}^{2*}\cdots
\cheu_{pq}^{2n*}\\
&=& b_p^{2n}\bar{w}_{2n}\alpha_p(\bar{w}_{2n}^*)
\alpha_p(b_q^{2n}\bar{w}_{2n}\alpha_q(\bar{w}_{2n}^*))
c^\alpha(p,q) \alpha_{pq}(\bar{w}_{2n})\bar{w}_{2n}^*b_{pq}^{2n*} \\
&=& b_p^{2n}\bar{w}_{2n}\alpha_p(\bar{w}_{2n}^*
b_q^{2n}\bar{w}_{2n})\alpha_p\alpha_q(\bar{w}_{2n}^*)
c^\alpha(p,q) \alpha_{pq}(\bar{w}_{2n})\bar{w}_{2n}^*b_{pq}^{2n*} \\
&=& b_p^{2n}{\theta}_{2n}\circs\alpha_p\circs {\theta}_{2n}^{-1}(
b_q^{2n})\bar{w}_{2n}c^\alpha(p,q)\alpha_{pq}\circs \sigma_{m(p,q)}(\bar{w}_{2n}^*)
\alpha_{pq}(\bar{w}_{2n})\bar{w}_{2n}^*b_{pq}^{2n*} \\
&=& b_p^{2n}{\theta}_{2n}\circs\alpha_p\circs {\theta}_{2n}^{-1}(
b_q^{2n}){\theta}_{2n}(c^\alpha(p,q))\bar{w}_{2n}\alpha_{pq}\left(
\sigma_{m(p,q)}(\bar{w}_{2n}^*)
\bar{w}_{2n}\right)\bar{w}_{2n}^*b_{pq}^{2n*} \\
&=& b_p^{2n}{\theta}_{2n}\circs\alpha_p\circs {\theta}_{2n}^{-1}(
b_q^{2n}){\theta}_{2n}(c^\alpha(p,q)){\theta}_{2n}\circs \alpha_{pq}\circs 
{\theta}_{2n}^{-1}
\left(\bar{w}_{2n}\sigma_{m(p,q)}(\bar{w}_{2n}^*)
\right)b_{pq}^{2n*} \\
&=&b_p^{2n}{\theta}_{2n}\circs \alpha_p\circs {\theta}_{2n}^{-1}(b_q^{2n})
{\theta}_{2n}(c^\alpha(p,q))
{\theta}_{2n}\circs \alpha_{pq}\circs {\theta}_{2n}^{-1}(
{\theta}_{2n}(\tilu_{m(p,q)})\tilu_{m(p,q)}^*) b_{pq}^{2n*}.
\end{eqnarray*}
Since the canonical extension is
 continuous in  the $u$-topology,  
 $\theta_{2n}(\tilu_{m(p.q)})\tilu_{m(p,q)}^*$ 
 converges to  $\bar{\theta}_0(\tilu_{m(p,q)})\tilu_{m(p,q)}^*$. 
Similar results holds for $c^{2n-1}(p,q)$.
Then by $(n.7)$, we have 
\begin{eqnarray*}
 \lefteqn{\hat{u}_p^0\bar{\theta}_0\circs\alpha_p\circs\bar{\theta}_0^{-1}(\hat{u}_q^0)
\bar{\theta}_0(c^{\alpha}(p,q))
\bar{\theta}_0\circs\alpha_{pq}\circs \bar{\theta}_0^{-1}(
\bar{\theta}_0(\tilu_{m(p,q)})\tilu_{m(p,q)}^*) \hat{u}_{pq}^{0*}} \\ &&=
\hat{u}_p^1\bar{\theta}_1\circs \beta_p\circs \bar{\theta}_1^{-1}(\hat{u}_q^1)
\bar{\theta}_1\circs \beta_{pq}\circs\bar{\theta}_1^{-1}(
\bar{\theta}_1(\tilu_{m(p,q)})\tilu_{m(p,q)}^*) \hat{u}_{pq}^{1*}
\end{eqnarray*}
and finished the proof of Theorem \ref{thm:quasiclass}. \hfill$\Box$
\bigskip

\noindent
\textbf{Proof of Theorem \ref{thm:class}.}
Set $\bar{\theta}_0(\tilu_n) =\tilu_n^\alpha$, $\bar{\theta}_1(\tilu_n)
=\tilu_n^\beta$, 
$\hat{u}_{\secp }=\hat{u}_p^{1*}\hat{u}_p^0$, $p\in Q$, 
and 
$\hat{u}_n=\tilu_n^\beta\tilu_n^{\alpha*} $, $n\in H$.

We replace 
$(\alpha_g, v^\alpha(g,h))$ and $(\beta_g,1 )$ with 
$(\bar{\theta}_0\circs\alpha_g\circs\bar{\theta}_0^{-1},
\bar{\theta}_0(v^\alpha(g,h)) )$ and
$(\bar{\theta}_1\circs\beta_g\circs\bar{\theta}_1^{-1},1)$ respectively. 
Note that we do not have $\alpha_n=\beta_n$, $n\in H$, after this
replacement, and in fact   
we have $\beta_n=\Ad \hat{u}_n\circs \alpha_n$, 
$\tal_n=\Ad \tilu_n^\alpha$ and $\tbe_n=\Ad \tilu_n^\beta$.

Summarizing results in Theorem \ref{thm:quasiclass}, 
we have the following. 
$$\beta_{\secp}=\Ad \hat{u}_{\secp}\circs\alpha_{\secp}, \,\,\,
\hat{u}_{\secp}{\alpha}_{\secp}(\hat{u}_{\sq}){c}^\alpha(p,q)
{\alpha}_{\secpq}(\hat{u}^*_{m(p,q)})\hat{u}_{\secpq}^*=1,\,\, p,q\in Q.$$  

For $r\in Q $ and $n\in H$, we define $\hat{u}_{n \sr}
=\hat{u}_{n}{\alpha}_{n}(\hat{u}_{\sr}){v}^\alpha(n,\sr)$.  
Then $\beta_g=\Ad \hat{u}_g\circs\alpha_g$ holds for every $g\in G$. 
Since 
\begin{eqnarray*}
\hat{u}_{ n }{\alpha}_{n }(\hat{u}_{\sr})
\delv^\alpha(\sr,n)^*
{\alpha}_{\sr}(\hat{u}_{\sr^{-1}n\sr}^*) 
\hat{u}_{\sr}^* 
&=& 
\hat{u}_{ n }\tilu_{ n }^\alpha
\hat{u}_{\sr} \tilu_{n }^{\alpha*}
\delv^\alpha(\sr,n)^*
{\alpha}_{\sr}(\hat{u}_{\sr^{-1}n\sr}^*)
\hat{u}_{\sr}^* \\
&=&
\hat{u}_{ n}\tilu_{ n }^\alpha
\tbe_{\sr}\tal_{\sr}^{-1}\left(\tilu_{n }^{\alpha*}
\delv^\alpha(\sr,n)^*{\alpha}_{\sr}(\hat{u}_{\sr^{-1}n\sr}^*)\right)
\\
&=&\lambda(\sr, n)
\hat{u}_{n }\tilu_{n }^\alpha
\tbe_{\sr}\left(\tilu_{\sr^{-1}n\sr }^{\alpha*}
\hat{u}_{\sr^{-1}n\sr}^*
\right) \\ 
&=&\lambda(\sr, n)
\tilu_{n }^\beta
\tbe_{\sr}(\tilu_{\sr^{-1}n\sr }^{\beta*}) \\
&=&1, 
\end{eqnarray*}
we have   
$\hat{u}_{n } 
{\alpha}_{n}(\hat{u}_{\sr})
{v}^\alpha(n,\sr)=
\hat{u}_{\sr}{\alpha}_{\sr}(\hat{u}_{\sr^{-1} n \sr}){v}^\alpha(\sr,\sr^{-1}n\sr).$
Compare this result with Lemma \ref{lem:cocycle0}.
Then in a similar way as in the proof of Lemma \ref{lem:cocycle} and 
Lemma \ref{lem:2cocycle}, we can 
show $\hat{u}_{g}{\alpha}_g(\hat{u}_h)v^\alpha(g,h)\hat{u}_{g h}^*=1$ for
$g,h\in G$. Thus $\alpha$ and $\beta$ are strongly cocycle conjugate.
\hfill$\Box$

\section{Model actions}\label{sec:model}
Construction of model actions with given invariant is presented in
\cite{Su-Tak-RIMS}, \cite{KtST} by using 
groupoid theory.  
However essential point of use of groupoid theory is
to construct a right inverse for the Connes-Takesaki module map. Thus
it may be possible to construct model actions without groupoid theory
once one admits the existence of the right inverse of the module map \cite{Su-Tak-fields}. 
In this section, 
we present the construction of
model actions along this observation.

Let $\varphi$ be a dominant weight, and 
$M=M_\varphi\rtimes_\theta
\mathbb{R}$ be a continuous decomposition, 
and  $u(s)$ the implementing unitary for $\theta_t$. 
Let $G$ a discrete amenable group, and $\alpha$
an action of $G$ on $M$. 
We quickly review how to describe $\mathrm{Inv}(\alpha)$ in terms of
$M_\varphi$. 

We may assume that  
$\varphi\circs\alpha=\varphi$, and $\alpha_g(u(s))=u(s)$ by cocycle
perturbation \cite{Su-Tak-act}. In this case, $M_\varphi$ is invariant under $\alpha_g$,
and we denote by $\alpha'_g$ the restriction on $M_\varphi$. 
Then $\mathrm{Inv}(\alpha)$ is obtained as follows.
A Connes-Takesaki module is given by
$\mathrm{mod}(\alpha_g)=\alpha'_g|_{Z(M_\varphi)}$. 
Let $H:=\{g\in G\mid \alpha'_g\in \Int(M_\varphi) \}$, which is a
centrally trivial part of $\alpha$. 
Fix $v_n\in U(M_\varphi)$ with $\alpha_n'=\Ad v_n$. Then we get 
a characteristic cocycle 
$(\lambda,\mu)\in Z(G,H, U(Z(M_\varphi)))$ and $c_t(n)\in
Z^1(\mathbb{R}, U(Z(M_\varphi)))$ as follows.
$$\alpha'_g(v_{g^{-1}ng})=\lambda(g,n)v_n,\,\, v_mv_n=\mu(m,n)v_{mn},\,\,
\theta_t(v_n)=c_t(n)v_n.$$ 
In this case, $\alpha_n=\Ad v_n\circs\sigma^\varphi_{c(n)}$ holds, 
where $\sigma^\varphi_{c(n)}$ is an extended modular automorphism \cite{CT}.

Conversely, for a normal subgroup $H\subset G$, a homomorphism $\beta:
g\in G\rightarrow 
\beta_g\in \Aut_\theta(Z(M_\varphi))$,
$(\lambda,\mu)$, and $c_t(n)$, we will construct a model action $\gamma_g$
with $\mathrm{Inv}(\gamma)=(H, \beta_g, [\lambda,\mu,c])$.

By \cite[Corollary 1.3]{Su-Tak-fields}, the exact sequence
$$ 1 \longrightarrow\overline{\Int}(M)\longrightarrow \Aut(M)\overset{\md}{\longrightarrow}
\Aut_\theta(Z(M_\varphi))\rightarrow 1$$ is split. 
By regarding $\beta$ as a
faithful homomorphism from $G/\mathrm{Ker}(\beta)$ into $\Aut_\theta(Z(M_\varphi))$, we lift
$\beta$ as an action of $G$ on $M$ by the above splitting exact
sequence. We may assume $\varphi\circs\beta_g=\varphi$,
$\beta_g(u(s))=u(s)$. 

Let $\partial(u)_t:=u\theta_t(u^*)\in Z^1(\mathbb{R},U(Z(M_\varphi)))$
for $u\in U(Z(M_\varphi))$. 
Since 
$$c_t(m)c_t(n)=c_t(mn)\partial(\mu(m,n)^*)_t,\,\,
\beta_g(c_t(g^{-1}ng))=\partial(\lambda(g,n)^*)_tc_t(n),$$ 
we have 
$$\sigma_{c(m)^*}^\varphi\circs\sigma_{c(n)^*}^\varphi=\sigma^\varphi_{c(m)^*c(n)^*}=
\sigma^\varphi_{\partial(\mu)c(mn)^*}=\Ad \mu(m,n)\circs
\sigma^\varphi_{c(mn)^*},
$$ and  $$\beta_g\circs \sigma^\varphi_{c(g^{-1}ng)^*}\circs \beta_g^{-1}=
\sigma^\varphi_{\beta_g(c(g^{-1}ng)^*)}=\sigma^\varphi_{\partial(\lambda(g,n))c(n)^*}=
\Ad \lambda(g,n)\circs
\sigma^\varphi_{c(n)^*}.
$$
In particular, 
$(\sigma_{c(n)^*}^\varphi,\mu(m,n))$ is a cocycle crossed action of $H$. 

Let $R_0$ be the injective factor of type II$_1$ with a tracial state
$\tau$, and $\alpha^{(0)}_g$ a free action of $G$ on $R_0$. Let
$\alpha_m:=\sigma_{c(m)^*}^\varphi\otimes 
\alpha_m^{(0)}$, and take a twisted crossed product 
$N=(M\otimes R_0)\rtimes_{\alpha,\mu\otimes 1}H$. Let $v_n$ be the
implementing unitary for $\alpha$.
Since $\mathrm{Inv}(\alpha)$ is
trivial, $N$ and $M$ have a common flow of weights by \cite{Kw-Tak},
\cite{Se-flow}. Hence $N\cong M$ by the classification theorem of
injective factors.  Let $\psi:=(\varphi\otimes \tau)\circs E$, where $E$
is the canonical conditional expectation on $M$. Then $\psi$ is
dominant, and $N_\psi=(M_\varphi\otimes R_0)\otimes_{\alpha,\mu\otimes 1}H$.

Let $\gamma_g=\beta_g\otimes \alpha_g^{(0)}$. Since
\begin{eqnarray*}
&&\gamma_g\circs \alpha_{g^{-1}ng}\circs \gamma_g^{-1}=\Ad
(\lambda(g,n)\otimes 1)\circs \alpha_n, \\
&&\lambda(g,m)\lambda(g,n)\mu(m,n)\lambda(g,mn)^*=\beta_g(\mu(g^{-1}mg,g^{-1}ng)), \\
&&\lambda(gh,n)=\beta_g(\lambda(h,g^{-1}ng))\lambda(g,n), 
\end{eqnarray*}
we can extend
$\gamma_g$ to an action on $N$ by 
$\gamma_g(v_{g^{-1}ng})=(\lambda(g,n)\otimes 1)v_n$. (See Appendix \ref{sec:ext} on 
the existence of such extension.) 

By the definition of $\gamma_g$, it is trivial that
$\mathrm{mod}(\gamma_g)=\beta_g$. If $m\in H$, then 
$\gamma'_m=\gamma_m|_{N_\psi}$ is given by $\Ad v_m$. Indeed if $x\in
M_\varphi\otimes 
R_0$, then $$\gamma_m'(x)=\id\otimes
\alpha^{0}_m(x)=\sigma_{c(m)^*}^\varphi\otimes \alpha_m^{(0)}(x)=\Ad v_m(x).$$
Note $\beta_m=\id$  and 
$\sigma^\varphi_{c(m)^*}|_{M_\varphi}=\id $ for $m\in H$.
By $\lambda(m,n)=\mu(m,m^{-1}nm)\mu(n,m)^*$,  $m,n\in H$, we have
$$\gamma'_m(v_{m^{-1}nm})=\lambda(m,n)v_{n}=
\mu(m,m^{-1}nm)\mu(n,m)^*v_n=
v_mv_{m^{-1}nm}v_{m}^*.$$ 
Due to the freeness of $\alpha^{(0)}_g$, 
$\gamma'_g\in \Int(N_\psi)$ if and only if $g\in H$. 
By the definition of $N$ and $\gamma_g$, we have
$\gamma_g(v_{g^{-1}ng})=\lambda(g,n)v_n$ and 
$v_mv_n=\mu(m,n)v_{mn}$. 
We can verify $\theta_t\otimes\id(v_n)=(c_t(n)\otimes
1)v_n$ as follows.
\begin{eqnarray*}
\theta_t\otimes \id(v_n)&=&(u(t)\otimes 1)v_n(u(t)\otimes 1)^*=
(u(t)\otimes 1)\alpha_{n}(u(-t)\otimes 1)v_n \\
&=& 
(u(t)\sigma_{c(n)^*}^\varphi(u(-t))\otimes 1)v_n =
(u(t)c_{-t}(n)^*u(t)\otimes 1)v_n \\ &=&(\theta_t(c_{-t}(n)^*)\otimes 1)v_n
=(c_t(n)\otimes 1)v_n.
\end{eqnarray*}

The above argument shows that $\gamma_g$ realizes the given
invariant. 

Here we treat
only type III factors, however
this construction is valid for the type II case.
If we use results in \cite{FT}, we can generalize the above construction
for arbitrary faithful normal semifinite weights.

\section{Group actions on subfactors}\label{sec:subfactor}
In this section, we see that we can apply our previous argument for
group actions on subfactors. We briefly recall basic notations for group
actions of subfactors

Let $N\subset M$ be a strongly amenable subfactor of type II$_1$ in the
sense of \cite{Po-amen}, and $N\subset M\subset M_1\subset M_2\subset
\cdots$ the Jones tower. 
For $\alpha\in\Aut(M,N)$, 
$\Phi(\alpha)=\{\alpha|_{M'\cap M_k}\}_{k=0}^\infty$ denotes the Loi invariant \cite{Loi-auto}. 
Let $\Cnt_{r}(M,N)$ be a set of all non-strongly outer automorphisms \cite{CK}, and  
$\chi_a(M,N)=\left(\mathrm{Ker}(\Phi)\cap \Cnt_{r}(M,N)\right)/\Int(M,N)$ the algebraic
$\chi$-group \cite{Go-chi}. Since $N\subset M$ is strongly amenable, 
$\mathrm{Ker}(\Phi)=\overline{\Int}(M,N)$ and $\Cnt(M,N)=\Cnt_{r}(M,N)$
hold. (See \cite{Loi-auto}, \cite{M-III1}, \cite{Po-act}.)
Take $\alpha\in \mathrm{Ker}(\Phi)$ and $\sigma\in \Cnt_{r}(M,N)$. Let $0\ne a\in
M_k$ be an element such that $\sigma(x)a=ax$ holds for all $x\in
M$. Then there exists a unitary $u(\alpha,\sigma)\in N$ such that
$\alpha(a)=u(\alpha,\sigma)a$ \cite{M-apprIJM}, which 
does not depend on $a$. This $u(\alpha,\sigma)$
satisfies $\alpha\circs \sigma\circs \alpha^{-1}=\Ad
u(\alpha,\sigma)\circs \sigma$. See \cite{M-apprIJM} for more properties of
$u(\alpha,\sigma)$.

We assume the following. \\
$(1)$ The normalizer groups for $N\subset M$ and $M\subset M_1$ are trivial. \\
$(2)$ There exists a lifting
$\sigma:\chi_a(M,N)\rightarrow \Aut(M,N)$. \\
Typical examples of such subfactors are the Jones subfactors with principal
graph $A_{2k+1}$, $k\geq 2$, \cite{J-ind}. 

Let $(\alpha,v^\alpha(g,h))$ be a cocycle crossed action of $G$ 
on $N\subset M$ with trivial Loi invariant.
(This condition corresponds to the triviality of Connes-Takesaki modules
in the previous sections.) 
Let $H=\{h\in G\mid \alpha_g\in \Cnt(M,N)\}$. The $\nu$-invariant
$\nu(n)$ is given by $\nu(n)=[\alpha_n]\in \chi_a(M,N)$. Note that 
$\nu(gng^{-1})=\nu(n)$ holds by $\Phi(\alpha_g)=\id$.
Fix $u_n^\alpha \in N$ with $\alpha_n=\Ad u_n^\alpha\circs\sigma_{\nu(n)}$.
Then the characteristic invariant $[\lambda,\mu]$ for $\alpha$ is given
as follows. 
$$\alpha_{g}(u_{g^{-1}n
g}^\alpha)u(\alpha_g,\sigma_{\nu(h)})=\lambda(g,n)\delta v^\alpha(g,n)u_n^\alpha,\,\,\,
u_m^\alpha\sigma_{\nu(m)}(u_n^\alpha)=\mu(m,n)v^\alpha(m,n)u_{m n}^\alpha.$$
We remark that $[\lambda,\mu]$ may be different from usual
characteristic invariant. Let $\kappa(k,l):=u(\sigma_l,\sigma_k)^*$ be
the $\kappa$-invariant for $N\subset M$. (This notion comes from \cite{J1}.)
The only difference is the
following relation.
$$\lambda(m,n)=\mu(m,m^{-1}nm)\overline{\mu(n,m)\kappa(\nu(n),\nu(m))},\,\,
m,n\in H.$$

In \cite{M-apprIJM}
 and \cite{M-notes}, we show $\mathrm{Inv}(\alpha)=(H, [\lambda,\mu],\nu)$ is a complete
 cocycle conjugacy invariant for approximately
inner actions of discrete
amenable groups on subfactors with conditions (1) and (2)
under some restrictions, e.g, the triviality of the $\kappa$-invariant.
 However if we modify the argument in
the previous sections in a suitable way, we can get rid of these
restrictions.

Let $(\beta, v^\beta(m,n))$ be another cocycle crossed action of $G$
with $\mathrm{Inv}(\beta)=\mathrm{Inv}(\alpha)$. We choose $u^\beta_m$
which satisfies the same relation for $(\lambda,\mu)$.
If we put $w_n:=u^a_mu^{\beta*}_m$ for $m\in H$, then $\alpha_n= \Ad
w_n\beta_n$ and $w_m\alpha_m(w_n)v^\alpha(m,n)w_{mn}^*=v^\beta(m,n)$
holds. Hence the same conclusion in Lemma \ref{lem:innervanish} holds. 
As in the single factor case, we may assume that 
$\alpha_n=\beta_n$ is a genuine action of $H$. 

If  $0\ne a\in M_k$ satisfies $\sigma_{\nu(n)}(x)a=ax$ for all $x\in M$, then
$a_n:=u_n^\alpha a$ satisfies $\alpha_n(x)a_n=a_nx$. Moreover, by the
definition of $\lambda(g,n)$ and $\alpha_g(a)=
u(\alpha_g,\sigma_{\nu(n)})a$, we have
$\alpha_{g}(a_{g^{-1}ng})=\lambda(g,n)\delta v^\alpha(g,n)a_n$. 
By using these facts, we can show a similar result 
in Lemma \ref{lem:cocycle0} as follows.

\begin{lem}\label{lem:cocyclesub} 
For any  $U_{g}=(v^\nu_g)\in  U(N^\omega)$ with 
 $\lim\limits_{\nu\rightarrow \omega}\alpha_{g} =\Ad v^\nu_{g} \circs\beta_{g}$. 
Then we have 
$U_{g}^*\delv^{\alpha}(g, n)\sigma_n(U_{g})=
\delv^{\beta}(g, n)$.
\end{lem}
\textbf{Proof.} Let $0\ne a_n\in M_k$ as above. 
Then 
\begin{eqnarray*}
 \lim_{\nu\rightarrow \omega }u_g^{\nu*}\delta v^\alpha(g,h)\alpha_n(u_g^\nu)a_n&=&
 \lim_{\nu\rightarrow \omega }u_g^{\nu*}\delta v^\alpha(g,h)a_nu_g^\nu 
=\beta_g\alpha_g^{-1}(\delta v^\alpha(g,h)a_n) \\
&=&\lambda(g,n)^*\beta_g(a_{g^{-1}ng}) 
=\delta v^\beta(g,n)a_n
\end{eqnarray*}
holds. 

Let $E$ be
the minimal conditional expectation from $M_k$ to $M$. Since $0\ne a_n a^*_n\in
M'\cap M_k$, $0\ne E(a_n a^*_n)\in \mathbb{C}$ follows, and hence we get the conclusion.
\hfill$\Box$
\smallskip

Thus we can repeat the same argument in the
previous sections, and classify approximately inner actions of discrete amenable 
groups on $N\subset M$. In particular, 
the higher obstruction
introduced in \cite{Kw-appr} and the $\nu$-invariant are complete outer
conjugacy invariants for automorphisms on Jones subfactors with
principal graph $A_{2k+1}$, $k\geq 2$. Note all automorphisms of 
Jones subfactors with principal graph $A_n$,
$n\geq 4$, are approximately inner. 

We close this section by explaining the construction of model actions.

Let $N\subset M$ be as above, and $R_0$ the injective
factor of type II$_1$. Let $\alpha_g^{(0)}$ be a  free action of $G$ on $R_0$.
Let $\alpha_n:=\sigma_{\nu(n)^{-1}}\otimes \alpha_n^{(0)} $, and 
$\tilde{\mu}(m,n):=\mu(m,n)\kappa(\nu(m),\nu(n))$. Take a twisted
crossed product $A\subset B:=(N\otimes R_0\subset M\otimes R_0)\rtimes_{\alpha,
\tilde{\mu}}H$. Since $\mathrm{Inv}(\alpha)$ is trivial, the standard invariant
of $A\subset B$ and that of $N\subset M$ are coincide, thus they are
isomorphic by Popa's classification theorem \cite{Po-amen}. Let $v_n$ be
the implementing unitary, and 
define $\gamma_g\in \Aut(B,A)$ as follows.
$$\gamma_g(x)=\id\otimes \alpha_g^{(0)}(x),\,\,x\in M\otimes R_0,\,\, 
\gamma_g(v_{g^{-1}ng})=\lambda(g,n)v_n.$$

We extend $\sigma_k\otimes \id$ to $\tilde{\sigma}_k\in \Aut(B,A)$ by 
$\tilde{\sigma}_k(v_n)=\overline{\kappa(k,\nu(n))}v_n$.
We claim that $\tilde{\sigma}_k$ is non-strongly outer. 
Take $0\ne a\in M_k$ such that
$\sigma_k(x)a=ax$, $x\in M$. It is clear $\tilde{\sigma}_k(x)(a\otimes
1)=(a\otimes 1)x$ for $x\in M_k\otimes R_0$. By the definition of the
$\kappa$-invariant, 
$$\tilde{\sigma}_k(v_n)(a\otimes
1)=\overline{\kappa(k,\nu(n))}
v_n(a\otimes 1)
=(u(\sigma_{\nu(n)},k)
\sigma_{\nu(n)^{-1}}(a)\otimes
1)v_n=(a\otimes 1)v_n.
$$
Thus $\tilde{\sigma}_k$ is non-strongly outer. 

We next show $\Ad v_n\tilde{\sigma}_n=\gamma_n$, $n\in H$. For $x\in M\otimes R_0$,
$$\Ad v_n\circs \tilde{\sigma}_{\nu(n)}=\Ad v_n\circs \sigma_{\nu(n)}\otimes
\id (x)=\id\otimes \alpha_n^{0}(x)=\gamma_n(x).$$
If we use
$\lambda(m,n)=\mu(m,m^{-1}nm)\overline{\mu(n,m)\kappa(\nu(n),\nu(m))}$, 
\begin{eqnarray*}
 \Ad v_m\circs \tilde{\sigma}_{\nu(m)}(v_{m^{-1}nm})&=&
\overline{\kappa(\nu(m),\mu(n))}v_mv_{m^{-1}nm}v_m^*\\&=&
\overline{\kappa(\nu(m),\mu(m^{-1}nm))}\tilde{\mu}(m,m^{-1}nm)\overline{\tilde{\mu}(n,m)}
v_{n}\\ &=&\mu(m,m^{-1}nm)\overline{\mu(n,m)\kappa(\nu(n),\nu(m))}v_n=\lambda(m,n)v_n,
\end{eqnarray*}
and hence we have $\Ad v_n\tilde{\sigma}_{\nu(n)}=\gamma_n$. We can
easily see that $\gamma_g$ is strongly outer for $g\in G\backslash H$.
Since
$\gamma_g(a\otimes 1)=a\otimes 1$, $u(\gamma_g,\tilde{\sigma}_k)=1$. By
the definition of $\gamma$, it is trivial that $[\lambda,\mu]$ is a
characteristic invariant for $\gamma$.

\section{$G$-kernels, or outer actions.}\label{sec:outer}
Let $M$ be an injective factor, and $G$ a discrete amenable group.
Let $\alpha$ be an injective homomorphism from $G$ to $\mathrm{Out}(M)$, 
or equivalently $\alpha$ be a map from $G$ to $\Aut(M)$ such that  
$\alpha_g\circs \alpha_h\equiv \alpha_{gh}\mod(\Int(M))$, and 
$\alpha_g\not\in \Int(M)$, $g\ne e$.
Such $\alpha$ is called a $G$-kernel, 
 or a free  outer action in \cite{KtT-outerI}.
In this section, we briefly explain that the intertwining argument is
applicable for classification of outer actions. (We always assume freeness.)
Difference of our argument with that of \cite{KtT-outerI} is
that we do not have to use a resolution group, which depends on the
choice of a representative 3-cocycle \cite[Remark 2.15]{KtT-outerI}.

We first recall invariants for outer actions introduced in
\cite{KtT-outerI}. 
Let $$H_\alpha=\{g\in G\mid \alpha_g\in \Cnt_r(M)\}.$$ 
As in the usual group
action case, we get $c_t(n)\in Z^1(\mathbb{R},Z(U(\tM)))$, 
$n\in H_\alpha$, by $\theta_t(\tilu_n)=c_t(n)\tilu_n$. Thus 
we get $\nu(n)=[c_t(n)]\in 
H^1(\mathbb{R},Z(U(\tM)))$. 
Set $Q=G/H_\alpha$ and fix a section $\secp\in G$ for $p\in Q$.
Fix $w^\alpha(p,q)\in U(\tM)$ such that  
$\tal_{\secp}\tal_{\sq}=\Ad w^\alpha(p,q)\circs\tal_{\secpq}$, $w(e,q)=w(p,e)=1$. 
Then we get $d_1^\alpha(p,q,r), d^\alpha_2(s;q,r) \in Z(U(\tM))$ by 
 $$w^\alpha(p,q)w^\alpha(pq,r)=d_1^\alpha(p,q,r)\tal_{\secp}(w^\alpha(q,r))w^\alpha(p,qr),\,\, 
\theta_s(w^\alpha(q,r))=d^\alpha_2(s;q,r)w^\alpha(q,r).$$ 
The modular obstruction $\mathrm{\Obm}(\alpha)$ is defined as
$([d_1^\alpha(p,q,r)\tal_{\secp}(d^\alpha_2(s;q,r)^*)],\nu)$.

\begin{thm}\label{thm:outerclass}
 Let $M$ and $G$ be as above. Let $\alpha$ and $\beta$ be outer actions
 of $G$ on $M$. If $(H_\alpha,\md(\alpha),\Obm(\alpha))=(H_\beta,\md(\beta),\Obm(\beta))$, then 
 $\alpha_g\equiv\sigma\circs \beta_g \circs \sigma^{-1}\!\mod\Int(M)$ for some
 $\sigma\in \overline{\Int}(M)$.
\end{thm} 
In the rest of this section, we assume
$(H_\alpha,\md(\alpha), \Obm(\alpha))=(H_\beta,\md(\beta),
\Obm(\beta))$. 
Hence 
 we can assume
$\alpha_n=\beta_n$ for $n\in H=H_\alpha$ as in \S 
\ref{sec:cocycle}. We fix $\tilu_n\in U(\tM)$ with
$\tal_n=\Ad \tilu_n$. 

We take $v^\alpha(g,h)\in U(M)$ such that
$\alpha_g\circs \alpha_h=\Ad v^\alpha(g,h)\circs \alpha_{gh}$. 
Then we get
$[\gamma^\alpha]\in H^3(G,\mathbb{T})$ by
$$v^\alpha(g,h)v^\alpha(gh,k)=\gamma^\alpha(g,h,k)\alpha_g(v^\alpha(h,k))v^\alpha(g,hk).$$ 
By \cite[Lemma 2.11]{KtT-outerI}, $[\gamma^\alpha]$ is uniquely
determined by $\Obm(\alpha)$.
As in the previous sections, set
$$c^\alpha(p,q)=v^\alpha(\secp,\sq)v^{\alpha}\left(\secpq,m(p,q)\right)^*, \,\,
\delta v^\alpha(g,n)=v^\alpha(g, g^{-1}ng)v^{\alpha}(n,g)^*.$$

Once  we fix $v^\alpha(g,h)$ and $\gamma^\alpha(g,h,k)$,  we 
get $\lambda^\alpha(g,n), \mu^\alpha(m,n), c_{t}^\alpha(n)
\in U(Z(\tM))$ as in \S \ref{sec:pre}, that is, they are defined as
follows.  
$$\alpha_{g}(\tilu_{g^{-1}ng})=\lambda(g,n)\delta v^\alpha(g,n)\tilu_{n},\,\,
\tilu_m\tilu_n=\mu^\alpha(m,n)v^\alpha(m,m)\tilu_{mn},\,\, \theta_t(\tilu_n)=c_t(n)\tilu_n.$$
These unitaries enjoy the
following relations.
\begin{eqnarray*}
 \lambda^{\alpha}(g,n)^*\theta_t(\lambda^\alpha(g,n))&=&
c_t^{\alpha}(n)^*\tal_g(c^\alpha_{t}(g^{-1}ng)), \\
c^\alpha_t(m)c^\alpha_t(n)c^{\alpha}_t(mn)^*&=&\mu^\alpha(m,n)^*\theta_t(\mu^\alpha(m,n)), \\
\gamma^\alpha(l,m,n)\mu^\alpha(l,m)\mu^\alpha(lm,n)&=&\mu^\alpha(m,n)\mu^\alpha(l,mn), \\
\lambda^\alpha(gh,n)&=&\tal_g(\lambda^\alpha(h,g^{-1}ng))\lambda^\alpha(g,n) \\
&&\times \gamma^\alpha(g,g^{-1}ng, h)
\overline{\gamma^\alpha(n,g,h)\gamma^\alpha(g,h,h^{-1}g^{-1}ngh)}, \\
\lambda^\alpha(g,mn)\lambda^\alpha(g,m)^*\lambda^\alpha(g,n)^*&=&
\mu^\alpha(m,n)\alpha_g(\mu^\alpha(g^{-1}mg,g^{-1}ng)^*)\\ 
&&\times \overline{\gamma^\alpha(m,g,g^{-1}ng)} 
\gamma^\alpha(g,g^{-1}mg, g^{-1}ng)\gamma^\alpha(m,n,g), \\
\lambda^\alpha(m,n)&=&\mu^\alpha(m,m^{-1}nm)\mu^\alpha(n,m)^*.
\end{eqnarray*} 
Let $Z(G,H;\gamma^\alpha, c_t^\alpha)$ be a set of all $(\lambda,\mu)$ satisfying
the above relation for fixed $\gamma^\alpha$ and $c_t^\alpha$. 
Here we can see that for $(\lambda,\mu),(\lambda',\mu')\in Z(G,H;\gamma^\alpha,
c_t^\alpha)$, $(\lambda (\lambda')^*, \mu(\mu')^*)$ is a usual
$\mathbb{T}$-valued characteristic cocycle.

We first clarify the relation between $(\lambda^\alpha, \nu^\alpha)$ and
$\Obm(\alpha)$. 
\begin{lem}\label{lem:obs}
 For $p,q,r\in Q$, we have
\begin{eqnarray*}
&&
c^\alpha(p,q)\alpha_{\secpq}(\delta v^\alpha(\sr,m(p,q))^*)
c^\alpha(pq,r)\alpha_{\secpqr}(v^\alpha(m(pq,r),\sr^{-1}m(p,q)\sr)) \\ 
&&=
\hat{\gamma}^\alpha(p,q,r)
\alpha_{\secp}(c^\alpha(q,r))c^\alpha(p,qr)
\alpha_{\secpqr}(v^\alpha(m(p,qr),m(q,r))).
\end{eqnarray*}
for some $\hat{\gamma}^\alpha(p,q,r)\in \mathbb{T}$.
\end{lem}
\textbf{Proof.}
Compute $(\alpha_{\secp}\circs \alpha_{\sq})\circs \alpha_{\sr}=
\alpha_{\secp}\circs (\alpha_{\sq}\circs \alpha_{\sr})$. Then we have
\begin{eqnarray*}
\Ad&&\hspace{-20pt} \left(c^\alpha(p,q)\alpha_{\secpq}(\delta v^\alpha(\sr,m(p,q))^*)
c^\alpha(pq,r)\alpha_{\secpqr}(v^\alpha(m(pq,r),\sr^{-1}m(p,q)\sr))\right)\circs 
\alpha_{\secpqr}\circs \sigma_{m(p,q,r)} \\
&&=
\left(\alpha_{\secp}(c^\alpha(q,r))c^\alpha(p,qr)\right)\circs 
\alpha_{\secpqr}(v^\alpha(m(p,qr),m(q,r)))\circs \alpha_{\secpqr}\circs\sigma_{m(p,q,r)}
\end{eqnarray*}
and obtain the desired conclusion.
\hfill$\Box$

\medskip
\noindent
\textbf{Remark.} $\hat{\gamma}^\alpha$ depends only on
$\gamma^\alpha$. 
To see this, let $(\beta,v^\beta(g,h))$ be another outer action of $G$ with
$\gamma^\alpha=\gamma^\beta$. 
Let $j(x)=x^*$ be a conjugate linear isomorphism from $M$ to $\Mopp$. 
Define $\theta_g:=\alpha\otimes j\circs \beta_g\circs j^{-1}$,
$v^\theta(g,h
)=v^\alpha(g,h)\otimes j(v^\beta(g,h))$. Then $(\theta_{g},v^\theta(g,h))$ is
a cocycle crossed action of $G$ on $M\otimes \Mopp$ due to the conjugate
linearity of $j(x)$. By Lemma \ref{lem:K-kercocycle}, we have 
\begin{eqnarray*}
&&c^\theta(p,q)\theta_{\secpq}(\delta v^\theta(\sr,m(p,q))^*)
c^\theta(pq,r)\theta_{\secpqr}(v^\theta(m(pq,r),\sr^{-1}m(p,q)\sr)) \\ 
&&=
\theta_{\secp}(c^\theta(q,r))c^\theta(p,qr)
\theta_{\secpqr}(v^\theta(m(p,qr),m(q,r))).
\end{eqnarray*}
Again by the 
conjugate linearity of $j(x)$, we get
$\hat{\gamma}^\alpha(p,q,r)=\hat{\gamma}^\beta(p,q,r)$.  

We can also describe $\hat{\gamma}^\alpha$ concretely. 
With a bit of little effort, we can show that 
$\hat{\gamma}^\alpha(p,q,r)$ is given by 
\begin{eqnarray*}
\hat{\gamma}^\alpha(p,q,r)&=&
\overline{\gamma^\alpha(\secp, \sqr, m(q,r))\gamma^\alpha(\secpq,m(p,q),\sr)
\gamma^\alpha(\secpqr, m(pq,r), \sr^{-1}m(q,r)\sr)} \\
&&\times\gamma^\alpha(\secp,\sq,\sr)
\gamma^\alpha(\secpq,\sr, \sr^{-1}m(p,q)\sr)
\gamma^\alpha(\secpqr, m(p,qr), m(q,r)) .
\end{eqnarray*}
We can also show $\hat{\gamma}^\alpha=\hat{\gamma}^\beta$ by the
above formula.
However we never use it
in the rest of this paper. So we omit the
proof of the above formula. 

\medskip

We fix $w^\alpha(p,q)$ as 
$w^\alpha(p,q)=c^\alpha(p,q)\tal_{\secpq}(\tilu_{m(p,q)}).$ 

\begin{lem}\label{lem:obs2}
For $(\lambda^\alpha,\mu^\alpha)\in Z(G,H; \gamma^\alpha, c_t^\alpha)$ and
 $p,q,r\in Q$, 
define 
$$\delta_{\lambda,\mu}(p,q,r)=\alpha_{\secpq}\left(\lambda^\alpha(\sr,m(p,q))^*\right)
\alpha_{\secpqr}\left(\mu^\alpha(m(pq,r),\sr^{-1}m(p,q)\sr)\mu^\alpha(m(p,qr),m(q,r))^*\right).
$$
Then $d^\alpha_1(p,q,r)=\hat{\gamma}^\alpha(p,q,r)\delta_{\lambda^\alpha,\mu^\alpha}(p,q,r)$ holds.
\end{lem} 
\textbf{Proof.} 
We compute $w^\alpha(p,q)w^\alpha(pq,r)$ and
$\alpha_{\secp}(w^\alpha(q,r))w^\alpha(p,qr)$. On one hand, we have
\begin{eqnarray*}
\lefteqn{ w^\alpha(p,q)w^\alpha(pq,r)} \\
&=&
c^\alpha(p,q)\alpha_{\secpq}(\tilu_{m(p,q)})
c^\alpha(pq,r)\alpha_{\secpq}(\tilu_{m(pq,r)}) \\
&=&
c^\alpha(p,q)\alpha_{\secpq}\alpha_{\sr}\alpha_{\sr}^{-1}(\tilu_{m(p,q)})
c^\alpha(pq,r)\alpha_{\secpqr}(\tilu_{m(pq,r)}) \\
&=&
c^\alpha(p,q) c^\alpha(pq,r)
\alpha_{\secpqr}\sigma_{m(pq,r)}\alpha_{\sr}^{-1}(\tilu_{m(p,q)})
\alpha_{\secpqr}(\tilu_{m(pq,r)}) \\
&=&
c^\alpha(p,q) c^\alpha_{pq,r}
\alpha_{\secpqr}\left(\tilu_{m(pq,r)}\alpha_{\sr}^{-1}(\tilu_{m(p,q)})\right) \\
&=&
c^\alpha(p,q) c^\alpha(pq,r)
\alpha_{\secpqr}\left(\tilu_{m(pq,r)}
\alpha_{\sr}^{-1}(\lambda^\alpha(\sr,m(p,q))^*\delta v^\alpha(\sr,
m(p,q))^*)\tilu_{\sr^{-1} m(p,q)\sr} 
\right) \\
&=&
\alpha_{\secpq}(\lambda^\alpha(\sr,m(p,q))^*)
c^\alpha(p,q) c^\alpha(pq,r) \\
&&\times \alpha_{\secpqr}\left(\sigma_{m(pq,r)}\left(
\alpha_{\sr}^{-1}(\delta v^\alpha(\sr,
m(p,q))^*)\right)\tilu_{m(pq,r)}
\tilu_{\sr^{-1} m(p,q)\sr} 
\right) \\
&=&
\alpha_{\secpq}(\lambda^\alpha(\sr,m(p,q))^*)
c^\alpha(p,q)
\alpha_{\secpq}
(\delta v^\alpha(\sr,
m(p,q))^*)
c^\alpha(pq,r)
\alpha_{\secpqr}
\left(\tilu_{m(pq,r)}\tilu_{\sr^{-1} m(p,q)\sr} 
\right) \\
&=&
\alpha_{\secpq}(\lambda^\alpha(\sr,m(p,q))^*)
\alpha_{\secpqr}(\mu^\alpha(m(pq,r),\sr^{-1}m(p,q)\sr)) \\
&&\times c^\alpha(p,q)
\alpha_{\secpq}
(\delta v^\alpha(\sr,
m(p,q))^*)
c^\alpha(pq,r)
\alpha_{\secpqr}
\left(\tilu_{m(p,q,r)}v^\alpha(
m(pq,r),\sr^{-1}m(p,q)\sr)
\right).
\end{eqnarray*}
On the other hand, we have
\begin{eqnarray*}
\alpha_{\secp}(w^\alpha(q,r))w^\alpha(p,qr)&=&
\alpha_{\secp}\left(c^\alpha(q,r)\alpha_{\sqr}(\tilu_{m(q,r)})
\right)c^\alpha(p,qr)\alpha_{\secpqr}(\tilu_{m(p,qr)}) \\
&=&
\alpha_{\secp}\left(c^\alpha(q,r)\right)
\alpha_{\secp}\alpha_{\sqr}(\tilu_{m(q,r)})
c^\alpha(p,qr)\alpha_{\secpqr}(\tilu_{m(p,qr)}) \\
&=&
\alpha_{\secp}\left(c^\alpha(q,r)\right)c^\alpha(p,qr)
\alpha_{\secpqr}\sigma_{m(p,qr)}(\tilu_{m(q,r)})
\alpha_{\secpqr}(\tilu_{m(p,qr)}) \\
&=&
\alpha_{\secp}\left(c^\alpha(q,r)\right)c^\alpha(p,qr)
\alpha_{\secpqr}(\tilu_{m(p,qr)}\tilu_{m(q,r)})\\
&=&\alpha_{\secpqr}(\mu^\alpha(m(p,qr)m(q,r))) \\
&&\times
\alpha_{\secp}\left(c^\alpha(q,r)\right)c^\alpha(p,qr)
\alpha_{\secpqr}(v^\alpha(m(p,qr),m(q,r))
\tilu_{m(p,q,r)}
).
\end{eqnarray*}
By Lemma \ref{lem:obs}, we get the conclusion.
\hfill$\Box$

\medskip 
The map $\delta_{\lambda,\mu}$ is the generalization of  
the Huebschmann-Jones-Ratcliffe map.
The following is a part of the Huebschmann-Jones-Ratcliffe exact
sequence \cite[Theorem 2]{Hueb},
\cite[Proposition 4.2.5]{J-act}, \cite[Theorem 8.1]{Ra}. 

\begin{lem}\label{lem:HJR}
 Let $[\lambda,\mu]\in \Lambda (G,H,\mathbb{T})$ be a usual
 $\mathbb{T}$-valued characteristic invariant. Then
 $\delta_{\lambda,\mu}=1$ in $H^3(Q,\mathbb{T})$ if and only if 
$\lambda(g,n)=\xi(g,g^{-1}ng)\overline{\xi(n,g)}$ and
 $\mu(m,n)=\xi(m,n)$ for some 2-cocycle $\xi(g,h)\in Z^2(G,\mathbb{T})$.
\end{lem}

\begin{lem}\label{lem:outchara}
 Let $\alpha$ and $\beta$ be as in Theorem \ref{thm:outerclass}. Then we
 can choose 
$(\gamma^\alpha, \lambda^\alpha, \mu^\alpha,
c_t^\alpha)=(\gamma^\beta,\lambda^\beta, \mu^\beta,c^\beta_t)$ and
 $v^\alpha(m,n)=v^\beta(m,n)$, $m,n\in H$.
\end{lem}
\textbf{Proof.} Set $\partial(z)(p,q,r)=\alpha_{\secp}(z(q,r))z(p,qr)z(pq,r)^*z(p,q)^*$
for $z(p,q)\in U(Z(\tM))$. 
Let $w^\beta(p,q)=c^\beta(p,q)\beta_{\secpq}(\tilu_{m(p,q)})$. 
Since we have fixed $\tilu_n\in U(\tM)$ in such a way $\tal_n=\tbe_n=\Ad
\tilu_n$, 
we have $c_t(n)=c^\alpha_t(n)=c^\beta_t(n)$. Take $z(p,q)\in U(Z(\tM))$ such
that $d^\alpha_1(p,q,r)=d^\beta_1(p,q,r)\partial(z)$ and $d^\alpha_2(s;
q,r)=d^\beta_2(s; q,r)z(q,r)\theta_s(z(q,r)^*)$. 
Since 
$\md(\alpha)=\md(\beta)$, 
\begin{eqnarray*}
w^{\alpha}(q,r)^*\theta_s(w^\alpha(q,r))&=&
d^\alpha_2(s; q,r)=\md(\alpha_{\sqr})(c_s(m(q,r)))\\ &=&d^\beta_2(s; q,r)
=w^{\beta}(q,r)^*\theta_s(w^\beta(q,r)) 
\end{eqnarray*}
Hence $z(q,r)=\theta_s(z(q,r))$ holds for all  $s\in \mathbb{R}$. By the ergodicity
of $\theta_s$, $z(q,r)\in \mathbb{T}$. By replacing $v^\beta(\secp,\sq)$
with $z(p,q)v^\beta(\secp,\sq)$, we can assume 
$d^\alpha_1(p,q,r)=d^\beta_1(p,q,r)$. 

Take 
 $a(g,h)\in \mathbb{T}$ such that $\gamma^\alpha=\gamma^\beta \partial(a)$.
By replacing $v^\beta(g,h)$
with $a(g,h)v^\beta(g,h)$, we can assume 
$\gamma^\alpha(g,h,k)=\gamma^\beta(g,h,k)=:\gamma(g,h,k)$. 
After this replacement 
we have $d^\beta_1=d^\alpha_2 \partial(b)$, where 
$b(p,q)=a(\secp,\sq)\overline{a(\secpq,m(p,q))}$. 
In this stage, we have
$(\lambda^\alpha,\mu^\alpha),(\lambda^\beta,\mu^\beta)\in Z(G,H;\gamma,c_t)$.

By Lemma
\ref{lem:obs2} and the remark after Lemma \ref{lem:obs},
$\delta_{\lambda^\beta,\mu^\beta}=\delta_{\lambda^\alpha,\mu^\alpha}\partial(b)$.  
Set $\lambda=\lambda^\alpha\lambda^{\beta*}$ and
$\mu=\mu^\alpha\mu^{\beta*}$. Then 
$(\lambda,\mu)$ is a usual $\mathbb{C}$-valued characteristic invariant with
$\delta_{\lambda,\mu}=1 $ in $H^3(Q,\mathbb{T})$. 
Let $\xi(g,h)$ be as in  Lemma \ref{lem:HJR}. By replacing
$v^\beta(g,h)$ with $\xi(g,h)v^\beta(g,h)$, we have
$(\lambda^\alpha,\mu^\alpha)=(\lambda^\beta, \mu^\beta)$. Note $\gamma(g,h,k)$
remains to be unchanged since $\xi(g,h)$ is a 2-cocycle, i.e.,
$\partial(\xi)=1$. 
Finally, $v^\alpha(m,n)=v^\beta(m,n)$ follows from
the definition of $\mu^\alpha(m,n)$ and $\mu^\alpha=\mu^\beta$. 
\hfill$\Box$
\smallskip

Once we have established Lemma \ref{lem:outchara}, we can repeat the similar
argument in \S \ref{sec:cocycle}, \S \ref{sec:ultra} and apply the
intertwining argument. (Though we can not
take $v^\alpha(m,n)=1$, $m,n\in H$, in this case, this  does not make
any problem.) Then we get $\alpha_{\secp}\equiv
\sigma\circs\beta_{\secp}\circs\sigma^{-1}\!\mod\Int(M)$ for some $\sigma\in
\overline{\Int}(M)$, and therefore $\alpha_g\equiv
\sigma\circs\beta_g\circs \sigma^{-1}\!\mod
\Int(M)$.

\appendix

\section{Proof of Proposition \ref{prop:2coho} and 
Theorem \ref{thm:Rohlin}}\label{sec:appen}

In this appendix, we present the proof of  Proposition \ref{prop:2coho}
and Theorem \ref{thm:Rohlin}.

In \cite[Proposition 7.4]{Ocn-act}, Ocneanu proved the
second cohomology vanishing theorem for semiliftable and \textit{strongly
free} actions of a discrete amenable group 
on $M_\omega$ by means of the Rohlin
theorem. However if we use the ultraproduct technique, we can present a
simple proof of the second cohomology vanishing theorem for semiliftable cocycle crossed
actions on an ultraproduct algebra without use of the Rohlin type theorem.

\begin{lem}\label{lem:perturb}
Let $R$ be a (not necessary separable) von Neumann algebra of type
 II$_1$ with a tracial state $\tau$,  
$Q$  a discrete amenable group, 
$(\alpha,u(g,h))$ a cocycle crossed action of $Q$ on $R$. Then for
 any finite subset $F\subset Q$ and $\epsilon>0$, there exist unitaries
 $\{w_g\}_{g\in Q}\subset R$ such that $|w_g\alpha_g(w_h)u(g,h)w_{gh}^*-1|_1<\epsilon$
 for $g,h\in F$. 
\end{lem}
\textbf{Proof.} Let $S$ be an $(F\cup F^2,\epsilon/2)$-invariant subset
of $Q$. For $g\in Q$, we 
define a bijection $\ell(g)$ on $S$ so that $\ell(g)h=gh$ if $gh\in S$.
Fix a system of matrix units $\{e_{g,h}\}_{g,h\in S}\subset R$ and set
$N=\{e_{g,h}\}'\cap R$. We identify $R$ with $N\otimes \{e_{g,h}\}''$.
Fix a
unitary $u_g$ such that $\Ad u_g\circs \alpha_g(e_{h,k})=e_{h,k}$. 
By replacing $\{\alpha,u(g,h)\}$ with $\{\Ad
u_g\circs \alpha_g,u_g\alpha_g(u_h)u(g,h)u_{gh}^* \}$, we may assume
$(\alpha, u(g,h))$
is of the form $(\alpha_g\otimes \id, u(g,h)\otimes 1)$ on $N\otimes
\{e_{g,h}\}''$. 
Set $w_{g}=\sum_{h\in S}u(g,h)^*\otimes e_{\ell(g)h,h}$.
For $g,h\in Q$, set $S_{g,h}=S\cap h^{-1}S\cap (gh)^{-1}S$. By the
choice of $S$, we have 
$|S\backslash S_{g,h}|\leq \epsilon|S|/2$, $g,h\in F$.

We have  
\begin{eqnarray*}
\lefteqn{ w_g\alpha_g(w_h)(u(g,h)\otimes 1)w_{gh}^*}\\ &=&
\left(\sum_{k\in S}u(g,k)^*\otimes e_{\ell(g)k,k}\right)
\left(\sum_{l\in S}\alpha_g(u(h,l)^*)\otimes e_{\ell(h)l,l}\right) \\
&&\times \left(u(g,h)\otimes 1\right)
\left(\sum_{m\in S}u(gh,m)\otimes e_{m,\ell(gh)m}\right) \\
&=&
\sum_{l\in S}u(g,\ell(h)l)^*
\alpha_g(u(h,l)^*)u(g,h)u(gh,l)\otimes e_{\ell(g)\ell(h)l,\ell(gh)l}\\
&=&\sum_{l\in S_{g,h}}1\otimes e_{ghl,ghl}+
\sum_{l\in S\backslash
S_{g,h}}u(g,\ell(h)l)^*\alpha_g(u(h,l)^*)u(g,h)u(gh,l)\otimes
 e_{\ell(g)\ell(h)l,\ell(gh)l}.
\end{eqnarray*}

Hence
$$
|w_g\alpha_g(w_h)u(g,h)w_{gh}^*-1|_1 \leq
\sum_{l\in S\backslash S_{g,h}}|1\otimes e_{l,l}|_1+
\sum_{l\in S\backslash S_{g,h}}|1\otimes e_{\ell(g)\ell(h)l,\ell(gh)l}|_1 
\leq \epsilon
$$
for $g,h\in F$. \hfill$\Box$

\smallskip
\noindent
\textbf{Proof of Proposition \ref{prop:2coho}.} 
Let $F_n$ be a finite subset of $Q$ with $F_n\subset
F_{n+1}$ and $\bigcup_n F_n =Q $. 
By Lemma \ref{lem:perturb}, there exist unitaries
$\{w_g^n\}$ such that
$|w_g^n\alpha_g(w_h^n)u(g,h)w_{gh}^{n*}-1|_1<1/n$ for $g,h\in F_n$.
Set $u^n(g,h)=w_g^n\alpha_g(w_h^n)u(g,h)w_{gh}^{n*}$. Then
$\lim_nu^n(g,h)=1$.
By the Index Selection Trick \cite[Lemma 5.5]{Ocn-act}, we get 
desired unitaries $\{w_g\}\subset M_\omega$. \hfill$\Box$  

\medskip

Next we present the proof of  Theorem \ref{thm:Rohlin},
which is simplification of that of  \cite[Theorem 6.1]{Ocn-act}.
Suppose  a countable subset $N\subset M_\omega$ is given. We may
assume $N$ is invariant under $\gamma$. 
Then $\gamma$ becomes a free proper action of $K$ on $N'\cap M_\omega$.

The following lemma due to Ocneanu \cite[Lemma 6.3]{Ocn-act} is out starting point. 
\begin{lem}\label{lem:Roh1}
 Let $\epsilon>0$ and $T\Subset K$ be such that $e\not\in T$. Then there
 exists a partition of unity $\{e_j\}_{j=0}^m\subset N'\cap M_\omega$ 
such that \\
$(1)$ $|e_0|_1\leq \epsilon$, \\
$(2)$ $e_i\gamma_g(e_i)=0$, $1\leq i\leq m$, $g\in T$.
\end{lem}

Take $F, S, \delta $ be as in Theorem \ref{thm:Rohlin}. 
Let  
$E=\{E_h\}_{h\in S}\subset N'\cap M_\omega$ be a set of mutually orthogonal projections. 
We set 
$$a_{g,E}=\sum_{h\in S\cap
g^{-1}S}|\gamma_g(E_h)-E_{gh}|_1, \,\,
b_E=\sum_{h\in S}|E_h|_1, \,\,
c_{g,E}=\sum_{h\in S\backslash g^{-1}S}|E_h|_1.$$

\begin{lem}\label{lem:order}
 Let $E=\{E_h\}_{h\in S}\subset N'\cap M_\omega$ be a set of mutually
 orthogonal projections with $[\gamma_g(E_h),E_k]=0$, $g\in F$, $h,k\in S$. 
If $b_E<1-\delta^{\frac{1}{2}}$, then there exists a 
set of mutually orthogonal projections $E'=\{E'_h\}_{h\in S}\subset
 N'\cap M_\omega$ such that 
\begin{enumerate}
\itemsep=0pt
\renewcommand{\labelenumi}{$(\arabic{enumi})$}
 \item $0<2^{-1}\delta^{\frac{1}{2}}\sum_{h\in S}
 |E_h'-E_h|_1\leq b_{E'}-b_E$,
 \item $a_{g,E'}-a_{g,E}\leq 2\delta^\frac{1}{2}(b_{E'}-b_E)$, $g\in F$,
 \item $c_{g,E'}-c_{g,E}\leq 2\delta^\frac{1}{2}(b_{E'}-b_E)$, $g\in F$,
 \item $[\gamma_g(E'_h),E_k']=0$, $g\in F, h,k\in S$.
\end{enumerate}
\end{lem}
\textbf{Proof.}
Fix $\epsilon>0$ with $b_E<(1-\delta^{\frac{1}{2}})(1-\epsilon)$. By
applying Lemma \ref{lem:Roh1} for $(FS)^{-1}FS\backslash \{e\}$, 
we get a partition of unity
$\{e_j\}_{j=0}^m$ in $N'\cap \{\gamma_g(E_h)\}_{g\in K,h\in S}'\cap
M_\omega$ such that $|e_0|_1<\epsilon$ and
$\gamma_g(e_j)\gamma_h(e_j)=0$, $1\leq j\leq m$, $g,h\in FS$, $g\ne h$.

Set $x=|S|^{-1}\sum_{g,h\in S}\gamma_g^{-1}(E_h)$. Then $|x|_1=\tau_\omega(x)=b_E$.
Here assume that  $|e_jx|_1\geq (1-\delta^{\frac{1}{2}})|e_j|_1$ for
$1\leq j\leq m$.
Then \begin{eqnarray*}
  b_E &=&|x|_1 
\geq  |(1-e_0)x|_1 
= \sum_{j=1}^m|e_jx|_1 
\geq \sum_{j=1}^m(1-\delta^{\frac{1}{2}})|e_j|_1 \\
&=&(1-\delta^{\frac{1}{2}})|1-e_0|_1 
\geq (1-\delta^{\frac{1}{2}})(1-\epsilon) 
>b_E
\end{eqnarray*}
holds by the choice of $\epsilon$. Hence 
we have $|e_jx|_1< (1-\delta^{\frac{1}{2}})|e_j|_1$ for some $j$.
Set $f=e_j$, $\tilde{f}=\sum_{g\in S}\gamma_g(f)$,
and $\rho=|\tilde{f}|_1=|S||f|_1$. 
Since 
$|fx|_1< (1-\delta^{\frac{1}{2}})|f|_1$,
\begin{eqnarray*}
 \left|\tilde{f}\sum_{h\in S}E_h\right|_1&=&\sum_{g\in S}\left|\gamma_g(f)\sum_{h\in S}
E_h\right|_1 
=\sum_{g\in S}\left|f\sum_{h\in S}\gamma_g^{-1}(E_h)\right|_1 \\
&=&|S||fx|_1 
< |S|(1-\delta^{\frac{1}{2}})|f|_1 
\leq (1-\delta^{\frac{1}{2}})\rho
\end{eqnarray*}
holds.

Set $f_h=\gamma_h(f)$, and 
$E'_h=E_h(1-\tilde{f})+\gamma_h(f)$. We show $E_h'\ne E_h$. If
we assume $E_h'=E_h$, then we get $\gamma_h(f)=E_h\tilde{f}$.  
However this is impossible because we have
$$(1-\delta^\frac{1}{2})\rho> \sum_{h\in
S}|\tilde{f}E_h|_1=\sum_{h\in S}|\gamma_h(f)|_1=\rho.$$

Since
$\{\gamma_{g}(f_h)\}$ are orthogonal projections for $g\in F, h\in S$,
commute with $\gamma_g(E_h)$ and $[\gamma_g(E_h), E_k]=0$,
 it follows that $[\gamma_g(E_h'),E_k']=0$
for $g\in F$, $h,k\in S$. 

We first verify the condition (1). Since 
$$
 \sum_{h\in S}|E_h'-E_h|_1 =\sum_{h\in S}|E_h\tilde{f}-f_h|_1 
\leq|\tilde{f}|_1+\sum_{h\in S}|f_h|_1 
\leq  2\rho $$
and 
\begin{eqnarray*}
 b_{E'}&=&\sum_{h\in S}|E_h'|_1 
= \sum_{h\in S}|E_h(1-\tilde{f})+f_h|_1 
= \sum_{h\in S}|E_h|_1+\sum_{h\in S}|f_h|_1-\sum_{h\in S}|E_h\tilde{f}|_1 \\
&\geq & b_E+\rho-(1-\delta^\frac{1}{2})\rho
=b_E+\delta^\frac{1}{2}\rho,
\end{eqnarray*}
we have
$$b_E'-b_E\geq \delta^\frac{1}{2}\rho\geq
\delta^{\frac{1}{2}}/2\sum_{h\in S}|E_h'-E_h|_1,$$
and we get the condition (1).

We next verify the condition (2). For $h\in S\cap g^{-1}S$, 
\begin{eqnarray*}
 |\gamma_g(E_h')-E_{gh}'|_1&=& 
 |\gamma_g(E_h(1-\tilde{f})+f_h)-E_{gh}(1-\tilde{f})-f_{gh}|_1 \\ 
&=& 
 |\gamma_g(E_h(1-\tilde{f}))-E_{gh}(1-\tilde{f})|_1 \\ 
&\leq &
 |\gamma_g(E_h)(1-\tilde{f})-E_{gh}(1-\tilde{f})|_1 
+ |\gamma_g(E_h)(\tilde{f}-\gamma_g(\tilde{f}))|_1 \\ 
&\leq &
 |\gamma_g(E_h)-E_{gh}|_1 
+ |\gamma_g(E_h)(\tilde{f}-\gamma_g(\tilde{f}))|_1  
\end{eqnarray*}
holds. Hence we have
\begin{eqnarray*}
 a_{g,E'}&=&\sum_{h\in S\cap g^{-1}S}
 |\gamma_g(E_h')-E_{gh}'|_1 \\
&\leq & 
\sum_{h\in S\cap g^{-1}S}
 |\gamma_g(E_h)-E_{gh}|_1 
+ |\gamma_g(E_h)(\tilde{f}-\gamma_g(\tilde{f}))|_1 \\ 
&\leq&a_{g,E}+
|\tilde{f}-\gamma_g(\tilde{f})|_1.  
\end{eqnarray*}
Here 
$$
 |\tilde{f}-\gamma_g(\tilde{f})|_1 =
\left|\sum_{h\in S}\gamma_h(f)-\sum_{h\in S}\gamma_g\gamma_h(f)\right|_1 
 \leq\sum_{h\in S \bigtriangleup gS}|\gamma_h(f)|_1 
\leq 2\delta|S||f|_1 
= 2\delta\rho
$$
holds. These inequalities yield
$$a_{g,E'}-a_{g,E}\leq 2\delta\rho\leq 2\delta^\frac{1}{2}(b_{E'}-b_E),
$$ and the condition (2) holds.

We finally verify the condition (3). Since we have
\begin{eqnarray*}
 c_{g,E'}&=& \sum_{h\in S\backslash g^{-1}S}|E_h'|_1 \\
&=&\sum_{h\in k\backslash g^{-1}S}|E_h(1-\tilde{f})|_1+
\sum_{h\in S\backslash g^{-1}S}|\gamma_h(f)|_1 \\
&\leq &c_{g,E}+|S\backslash g^{-1}S||f|_1 \\
&\leq & c_{g,E}+2\delta|S||f|_1 \\
&\leq & c_{g,E}+2\delta\rho \\
&\leq & c_{g,E}+2\delta^\frac{1}{2}(b_{E'}-b_E),
\end{eqnarray*}
we get the condition (3) \hfill$\Box$

\noindent
\textbf{Proof of Theorem \ref{thm:Rohlin}.}
Let $\mathcal{S}$ be a set of families of orthogonal projections
$E=\{E_h\}_{h\in S}$ in $N'\cap M_\omega$ satisfying the following
conditions.
\begin{eqnarray*}
 (1)&& a_{g,E}\leq 2\delta^\frac{1}{2}b_E, \quad g\in F, \\
 (2) &&c_{g,E}\leq 2\delta^\frac{1}{2}b_E, \quad g\in F, \\
 (3)&& [\gamma_g(E_h),E_k]=0, \quad g\in F,\,\, h,k\in S.
\end{eqnarray*}
Obviously $\mathcal{S}$ is not empty.
We  define an order $E\leq E'$ on $\mathcal{S}$ if $E=E'$ or $E$ and $E'$ satisfy 
the first three conditions in Lemma \ref{lem:order}. In the same way as
in \cite{Ocn-act}, 
it is shown that $\mathcal{S}$ is an inductive ordered set.
Let $\bar{E}=\{\bar{E}_h\}$ be a maximal element, and
assume $b_{\bar{E}}<1-\delta^\frac{1}{2}$.  By Lemma
\ref{lem:order}, there exists a partition of unity $\{E_h'\}$ satisfying  
the conditions in Lemma \ref{lem:order}. Then it is easy to see
$\{E_h'\}$ is in $\mathcal{S}$, and $\{E_h'\}$ is strictly larger than
$\bar{E}$.

Hence $b_{\bar{E}}\geq 1-\delta^\frac{1}{2}$ holds. Let
$E_0=1-\sum_{h\in S}\bar{E}_h$. Then $|E_0|_1\leq
\delta^{\frac{1}{2}}$.
Fix $h_0\in S$, and set $E_{h_0}=\bar{E}_{h_{0}}+E_0$,  and
$E_h=\bar{E}_h$ for $h\ne h_0$. It is obvious $[\gamma_g(E_h),E_k]=0$,
$g\in F$, $h,k\in S$. It is also easy to see 
$a_{g,E}\leq 4\delta^{\frac{1}{2}}$ and $c_{g,E}\leq
3\delta^\frac{1}{2}$. \hfill$\Box$

\section{Extension of automorphisms to a twisted crossed product algebra}\label{sec:ext}

Let $M$ be a von Neumann algebra, and $L^2(M)$ the standard
Hilbert space. 
Let $(\alpha_g,v(g,h))$ be a cocycle crossed action of $G$. 
Define  $\pi(x),\lambda_g\in B(L^2(M)\otimes \ell^2(G))$, $x\in M$, $g\in G$, as follows. 
(Do not confuse $\alpha_{g^{-1}}$ and $\alpha_g^{-1}$.)
$$\left(\pi(x)\xi\right)(g)=\alpha_{g^{-1}}(x)\xi(g),\,\,\left(\lambda_g\xi\right)(h)
=v(h^{-1},g)\xi(g^{-1}h).
$$
Then $\lambda_g\pi(x)\lambda_g^*=\pi(\alpha_g(x))$ and 
$\lambda_g\lambda_h=\pi(v(g,h))\lambda_{gh}$ hold. 
By definition, the twisted crossed
product von Neumann algebra is $M\rtimes_{\alpha,v}G=\pi(M)\vee\{\lambda_g\}$.

Although the following result may be well-known, we present a proof for
readers' convenience.

\begin{thm}\label{thm:ext}
$(1)$ Assume that $\theta\in \Aut(M)$ induces an
 automorphism on $G$ (denoted by the same symbol $\theta$) such that 
$$\theta\circs\alpha_{\theta^{-1}(g)}\circs \theta^{-1}=\Ad
 u_g^\theta\circs \alpha_g(x),\,\, 
u_g^\theta\alpha_g(u_h^\theta)v(g,h)u_{gh}^{\theta*}=\theta(v(\theta^{-1}(g),\theta^{-1}(h)))$$
for some unitary $u_g^\theta\in M$.
Then there exists a unique extension $\tilde{\theta}\in \Aut(M\rtimes_{\alpha, v}G)$ of $\theta$ 
such that 
 $\theta(\lambda_{\theta^{-1}(g)})=\pi(u_g^\theta)\lambda_g$. \\
$(2)$
 If we have
 $u_g^{\theta\sigma}=\theta(u_{\theta^{-1}(g)}^\sigma)u_g^\theta$
for  $\theta,\sigma\in \Aut(M)$,  then
$\tilde{\theta}\circs \tilde{\sigma}=\widetilde{\theta\circs \sigma}$ holds.
\end{thm}
\textbf{Proof.} 
(1) Denote $u_g^\theta$ by $u_g$ for simplicity.
Let $U_\theta\in B(L^2(M))$ be the standard implementing
unitary of $\theta$.
Define a unitary $W_\theta\in B(L^2(M)\otimes \ell^2(G))$ by 
$$\left(W_\theta \xi\right)
(g)=u_{g^{-1}}^*
U_\theta \xi(\theta^{-1}(g)),\,\,\, \xi(g)\in L^2(M)\otimes \ell^2(G).$$ 
Then $W_\theta^*$ is given by 
$$\left(W_\theta \xi\right)
(\theta^{-1}(g))=U_\theta^*u_{g^{-1}}\xi(g).$$ 
We show $\Ad W_\theta$ gives a desired action.  
At first, we verify $\Ad W_\theta (\pi(a))=\pi(\theta(a))$ as follows.
\begin{eqnarray*}
 \left(W_\theta\pi(a)W_\theta^*\xi\right)(g)&=&
u_{g^{-1}}^*
U_\theta\left(\pi(a)W_\theta^*\xi\right)(\theta^{-1}(g)) \\
&=&
u_{g^{-1}}^*
U_\theta\alpha_{\theta^{-1}(g^{-1})}(a)\left(W_\theta^*\xi\right)(\theta^{-1}(g)) \\
&=&
u_{g^{-1}}^*
U_\theta\alpha_{\theta^{-1}(g^{-1})}(a)
U_\theta^*u_{g^{-1}}\xi(g)  \\
&=&
u_{g^{-1}}^*
\theta\alpha_{\theta^{-1}(g^{-1})}(a)u_{g^{-1}}\xi(g) \\ 
&=&
\alpha_{g^{-1}}(\theta(a))
\xi(g)  \\
&=&
\left(\pi(\theta(a))\xi\right)(g).
\end{eqnarray*}

Next we show $\Ad W_\theta(\lambda_{\theta^{-1}(g)})=\pi(u_g)\lambda_g$ as follows. 
\begin{eqnarray*}
\left( W_\theta \lambda_{\theta^{-1}(g)} W_\theta^*\xi \right)(h)
&=&
u_{h^{-1}}^*
U_\theta\left(\lambda_g W_\theta^*\xi
			  \right)(\theta^{-1}(h))\\
&=&
u_{h^{-1}}^*U_\theta\alpha_{\theta^{-1}(h^{-1})}\left(v(\theta^{-1}(h^{-1}),
\theta^{-1}(g))\right)
\left(W_\theta^*\xi \right)(\theta^{-1}(g^{-1}h))\\
&=&
u_{h^{-1}}^*U_\theta\alpha_{\theta^{-1}(h^{-1})}\left(v(\theta^{-1}(h^{-1}),
\theta^{-1}(g))\right)
U_\theta^*u_{h^{-1}g}
\xi(g^{-1}h)\\
&=&
u_{h^{-1}}^*\theta\alpha_{\theta^{-1}(h^{-1})}\left(v(\theta^{-1}(h^{-1}),
\theta^{-1}(g))\right)u_{h^{-1}g}
\xi(g^{-1}h)\\
&=&
u_{h^{-1}}^*u_{h^{-1}}\alpha_{h^{-1}}(u_g)
v(h^{-1},g)
\xi(g^{-1}h)\\
&=&
\alpha_{h^{-1}}(u_g)
\left(\lambda_g\xi\right)(h)\\
&=&
\left(\pi(u_g)\lambda_g\xi\right)(h).
\end{eqnarray*}
This shows that $\tilde{\theta}=\Ad W_\theta$  preserves $M\rtimes_{\alpha, v}G$, and 
is a desired extension. The uniqueness of $\tilde{\theta}$ is obvious.\\
(2) By the assumption, $W_\theta W_\sigma=W_{\theta\sigma}$ holds.
\hfill$\Box$
\smallskip

\ifx\undefined\bysame
\newcommand{\bysame}{\leavevmode\hbox to3em{\hrulefill}\,}
\fi

\end{document}